\documentclass[3p]{elsarticle}
\usepackage{graphicx}
\usepackage{cancel}
\usepackage{epstopdf}
\usepackage[english]{babel}
\usepackage{natbib}
\usepackage{hyperref}
\usepackage{tikz}
\usepackage{fullpage}
\usepackage{microtype}
\usepackage{algorithm,algorithmic}
\usepackage{mcode}

\author{Sumedh M.~Joshi}
\address{Center for Applied Mathematics, \\ 657 Rhodes Hall, Cornell University, \\ Ithaca NY, 14850}

\author{Greg N.~Thomsen}
\address{Applied Research Laboratories \\ 10000 Burnet Road,  \\ Austin TX, 78758}

\author{Peter J.~Diamessis}
\address{School of Civil and Environmental Engineering, \\ 105 Hollister Hall, Cornell University, \\ Ithaca NY, 14850}

\cortext[cor1]{Please send correspondence to the first author at smj96@cornell.edu}

\usepackage{amsmath}   
\usepackage{amsthm}
\usepackage{amssymb} 
\usepackage{calc}
\usepackage{ifthen}
\usepackage{graphicx}
\usepackage{enumerate}
\usepackage{centernot}
\usepackage[mathscr]{eucal}

\renewcommand\and{\ensuremath{\;\;\;\;\;\;\text{and}\;\;\;\;\;\;} }

\newcommand\restr[2]{{
  \left.\kern-\nulldelimiterspace 
  #1 
  \vphantom{\big|} 
  \right|_{#2} 
}}

 \newcommand\R{\ensuremath{\bb{R}}}

 \newcommand\C{\ensuremath{\bb{C}}}

\newcommand\norm[1]{\ensuremath{\left\lvert \left\lvert  #1 \right\rvert \right\rvert }}

\newcommand\avg[2][1]{
	\ifthenelse{\equal{#1}{1}}{\ensuremath{\left \langle#2\right \rangle }}{}
	\ifthenelse{\equal{#1}{2}}{\ensuremath{\left \langle#2^2\right \rangle }}{}
	\ifthenelse{\equal{#1}{3}}{\ensuremath{\left \langle#2\right \rangle ^2}}{}}

\newcommand\derv[3][1]{
	\ifthenelse{\equal{#1}{1}}{\ensuremath{\frac{d#2}{d#3}}}{}
	\ifthenelse{\equal{#1}{2}}{\ensuremath{\frac{d^2#2}{d#3^2}}}{}}

\newcommand\prtl[3][1]{	\ifthenelse{\equal{#1}{1}}{\ensuremath{\frac{\partial#2}{\partial#3}}}{}
	\ifthenelse{\equal{#1}{2}}{\ensuremath{\frac{\partial^2#2}{\partial#3^2}}}{}}

\renewcommand\bf[1]{\ensuremath{\mathbf{#1}}}
\newcommand\bb[1]{\ensuremath{\mathbb{#1}}}

\renewcommand\eqref[1]{Eq. (\ref{#1})}

\theoremstyle{plain}   
\theoremstyle{definition}  
\numberwithin{exer}{subsection} 
\theoremstyle{plain}   \newtheorem{cla}{Claim}
\theoremstyle{plain}   
\theoremstyle{definition}   
\theoremstyle{definition}   
\theoremstyle{plain}   
\theoremstyle{plain}   
\theoremstyle{plain}    

\title{Deflation-accelerated preconditioning of the Poisson-Neumann Schur problem on long domains with a high-order discontinuous element-based collocation method}

\begin{document}

\begin{abstract}
A combination of block-Jacobi and deflation preconditioning is used to solve a high-order discontinuous collocation-based discretization of the Schur complement of the Poisson-Neumann system as arises in the operator splitting of the incompressible Navier-Stokes equations.
The preconditioners and deflation vectors are chosen to mitigate the effects of ill-conditioning due to highly-elongated domains and to achieve GMRES convergence independent of the size of the grid.
The ill-posedness of the Poisson-Neumann system manifests as an inconsistency of the Schur complement problem, but it is shown that this can be accounted for with appropriate projections out of the null space of the Schur matrix without affecting the accuracy of the solution.
The combined deflation/block-Jacobi preconditioning is compared with two-level non-overlapping additive Schwarz preconditioning of the Schur problem, and while both methods achieve convergence independent of the grid size, deflation is shown to require half as many GMRES iterations and $25\%$ less wall-clock time for a variety of grid sizes and domain aspect ratios.
The deflation methods shown to be effective for the two-dimensional Poisson-Neumann problem are extensible to the three-dimensional problem assuming a Fourier discretization in the third dimension.
A Fourier discretization results in a two-dimensional Helmholtz problem for each Fourier component that is solved using deflation/block-Jacobi preconditioning on its Schur complement.
Here again deflation is shown to be superior to two-level non-overlapping additive Schwarz preconditioning, requiring about half as many GMRES iterations and $15\%$ less time.
While the methods here are demonstrated on a spectral multidomain penalty method discretization, they are readily extensible to any discontinuous element-based discretization of an elliptic problem, and are particularly well-suited for high-order methods.
\end{abstract}


\maketitle

\section{Introduction }

When operator-splitting the incompressible Navier-Stokes equations in the context of computational fluid dynamics, the satisfaction of the divergence-free condition on the velocity requires the solution of a Poisson equation with Neumann boundary conditions \cite{Karniadakis1991}.  This pressure Poisson equation plays the role of enforcing incompressibility of the velocity field, and computing its accurate solution is the dominant computational expense in numerically solving the incompressible Navier-Stokes equations\cite{Escobar-Vargas2014}.  This paper describes a method for solving the  Poisson-Neumann problem\footnote{in which the real-valued function $u(x,z)$ represents the pressure.},
	\begin{align}
		\nabla^2 u &= f  \textrm{    on    } \Omega \nonumber \\
		n\cdot \nabla u &= g \textrm{    on   } \partial\Omega,
		\label{poisson}
	\end{align}
as it arises in the context of solving the incompressible Navier-Stokes equation on element-based discrete grids in environmental fluid dynamics.  Due to the fact that any constant function added to a solution $u$ is also a solution, this problem is ill-posed\cite{Pozrikidis2001}, and so a compatibility condition,
\begin{align}
	\int_\Omega f d\Omega = \int_{\partial \Omega} g dS
\end{align}
must be added to make the system solvable.

While the discrete Poisson problem in general has certainly garnered much attention in the literature, the large domains with large aspect ratios in environmental flows add unique challenges from a computational perspective.  These difficulties can be characterized as follows.
First, to avoid artificial numerical dispersion and dissipation that can pollute a low-order discretization\cite{Vitousek2011,Diamessis2005}, high-order discretizations are required for accurately computing dynamical properties of environmental flows with reasonable computational effort.  These high-order methods are more difficult to implement in practice, and they can result in ill-conditioned operator matrices that must be carefully constructed to avoid numerical error \cite{Costa2000}.
Secondly, due to the presence of gravity, environmental problems often have very different spatial and temporal scales in the vertical direction than in the horizontal direction.  Thus, computational grids tend to be highly anisotropic with elements frequently compressed in the vertical direction and greatly elongated in the horizontal direction \cite{Scotti2008,Santilli2015,Vitousek2014}.  From a numerical perspective, high-order discretizations and high-aspect ratio elements lead to ill-conditioned linear systems and increased computational cost in obtaining their solution.

One approach to addressing the particular difficulties in environmental flows is to exploit physical properties to tailor a numerical method.  Asymptotic expansions in the leptic (aspect) ratio of the grid are sometimes used to divide the Poisson operator into two coupled problems, one in the vertical (thin) direction, and another problem in the remaining dimensions \cite{Scotti2008}.  These methods have been successfully combined with geometric multi-grid \cite{Santilli2015} and Krylov methods \cite{Santilli2011} for obtaining the solution to the independent sub-problems derived in the leptic expansion.  Yet another physical approach is to formulate the Poisson problem in a coordinate system following the perturbations of the horizontal iso-contours of the flow field.  By making a small-slope approximation of the iso-contours, a simpler Poisson problem is solved.  This method has been shown to be effective in reproducing the dynamics of nonlinear internal wave propagation at the environmental scale \cite{Vitousek2014}, but cannot capture more localized multiscale phenomena such as turbulence and the nonlinear stages of evolution of a two-dimensional instability.

Apart from such physically-inspired approaches, ideas from domain decomposition can also be used to address the difficulty of solving high-order long-domain discretizations at the algebraic level.  While domain decomposition methods are general, they have properties that are beneficial for the ill-conditioned Poisson-Neumann problem on stretched grids.  Since most environmental flow problems tend to be too large for direct matrix factorization, domain decomposition methods focus on constructing good preconditioners for the iterative solvers.  Efforts have been made to alleviate the ill-conditioning of the element matrices due to aspect-ratio, as well as to yield convergence of iterative solvers independent of the size of the number of elements in the grid \cite{Fischer2005,Feng2001}.  Since in environmental problems domains have high aspect ratio and often contain many elements, domain decomposition is well-suited for this application.

In this work, two-dimensional domain decomposition is used to construct decoupled local problems and a Schur complement system that, while smaller and better numerically conditioned than the full problem, is still generally too large in practice to solve directly.  Thus the solution of the Schur complement problem is usually obtained iteratively.  The Schur complement approach to solving the elliptic equations has been implemented for high-order \cite{Couzy1995} spectral and collocation-based \cite{Bialecki2007} discretizations with iterative Krylov solvers.  As is common with Krylov methods, the Schur complement approach ubiquitously requires effective preconditioners to be practical.  A modern class of domain-decomposition based preconditioners are the multilevel additive Schwarz methods, and they have been shown to be very effective on the Schur complement problem of the Helmholtz \cite{Manna2004,Pasquetti2006} and the Poisson-Dirichlet equations \cite{Pavarino2000}.

Closely related to multi-level additive Schwarz methods, deflation methods accelerate the convergence of Krylov methods by attempting to identify and directly solve components of the residual associated with slowly-converging eigenvectors \cite{Vuik2006}.
Deflation methods have been shown to be superior to two-level additive Schwarz methods for fluid-flow problems with symmetric positive-definite discretizations \cite{Nabben2004}, but to the knowledge of the authors  have not been compared with additive Schwarz methods for the Schur complement problem, nor been applied to environmental flow problems.
However, deflation methods have been shown to be effective for the Poisson-Neumann problem \cite{Tang2009} and so it is reasonable to expect deflation to perform well for the Schur complement problem as well.  In this paper, it is demonstrated that on high-order, high-aspect ratio grids deflation augmented with block-Jacobi preconditioning achieves convergence of the Poisson-Neumann Schur system independent of both the aspect ratio and the number of elements in the long direction.  Furthermore it is shown that deflation/block-Jacobi preconditioning is  at least as effective as the two-level additive Schwarz method, requiring approximately half as many Krylov iterations as the two-level additive Schwarz preconditioning method.   While this is demonstrated on the Spectral Multidomain Penalty Method \cite{Diamessis2005} discretization, the ideas outlined herein may be applied to any high-order discontinuous element-based discretization of an elliptic problem, which encompasses a broad and important class of methods and partial differential equations.

The paper is organized as follows.  In Section 2 is introduced the discretization of Poisson-Neumann problem used here along with the construction of its Schur complement problem. In Section 3 is described the deflation/preconditioning method for solving the Schur complement system. In Section 4 are the results comparing the performance of deflation augmented preconditioning with two-level additive Schwarz preconditioning on a suite of test problems.  In Section 5 is shown an extension of this method to the three-dimensional problems.  Section 6 is a discussion of the broader applicability of the results, and Section 7 is a summary of the results and concluding remarks.

\section{Construction of the Schur problem}
\label{sec:schur}
   %
\subsection{Background on Schur complement methods}
Schur complement methods for partial differential equations arise in many contexts related to domain decomposition of element-based grids.  Either as preconditioners or solvers, domain decomposition methods have been used extensively in solving large sparse discretizations of partial differential equations\cite{Toselli2004},  by separating the problem into easily parallelizable local interior problems and a separate Schur complement problem.   When used as a preconditioning technique, so-called additive Schwarz methods solve the interior problems directly for use as a preconditioner to the full problem, often with sub-domain overlap to achieve condition numbers that scale as $C(1 + (H/\delta)^2)$ where $H$ is the subdomain size, $\delta$ the overlap, and $C$ a constant independent of both \cite{Brenner1996,Hesthaven,Fischer2005,Feng2001}.  With the addition of a coarse-level correction, this bound can be reduced to $C(1 + H/\delta)$ on quasi-regular grids.  When used as a solver, domain decomposition techniques solve for the interface unknowns first using the Schur complement system, whose condition number scales as the inverse of the sub-domain size $h$, $Ch^{-1}$, and can be significantly smaller in dimension than the full operator \cite{Brenner1999}.  In doing so, the Schur problem is sometimes solved iteratively via multi-grid\cite{Beuchler2002}, conjugate gradients\cite{Cros2002,Carvalho2001}, or GMRES\cite{Yamazaki2010a,Kanschat2003a}.  Coarse-grid corrected additive-Schwarz preconditioners on the Schur problem reduce the condition number to $C(1 + \log(H/h)^2)$\cite{Carvalho2001}.  Thus, whenever possible, it is preferable to solve the Schur complement problem with adequate local and global (coarse-grid) preconditioning as opposed to solving a full Poisson system. 

\subsection{The spectral multi-domain penalty method}
\label{smpm}
To discretize the two-dimensional Poisson-Neumann system a high-order discontinuous collocation-based variant of the spectral element method called the Spectral Multidomain Penalty Method (SMPM)\cite{Diamessis2005,Escobar-Vargas2014} is used.
In the SMPM each element is assumed to be smoothly and invertibly mapped from the unit square $[-1,1]\times[-1,1]$ and the element connectivity is logically cartesian (each element has a single neighbor in each of the North, South, East, and West directions).
Within each element lies a two-dimensional Gauss-Lobatto-Legendre (GLL) grid; denote as $n$ the number of GLL points per direction per element, and $m_x$ and $m_z$ the number of $x$ and $z$ elements in the grid\footnote{Here $z$ is the vertical direction as is convention in environmental fluid mechanics.}.  
On the GLL grid, a two-dimensional nodal Lagrange interpolant basis of polynomial order $n+1$ is constructed such that each basis element is one on one of the $n^2$ GLL points and zero on all of the others.
This basis is used for approximating grid functions and their derivatives which are calculated by way of spectral differentiation matrices \cite{Costa2000} which compute derivatives of nodally-represented functions by multiplying the nodal values by derivatives of the Lagrange interpolants themselves.
The SMPM is a discontinuous method and so $C^0/C^1$ inter-element continuity and boundary conditions are only weakly enforced.
While this method has been shown to be stable, spectrally accurate, and effective for under-resolved high-Reynolds number simulations of incompressible flows \cite{Diamessis2005,Zhou2013,Diamessis2006,Abdilghanie2013}, the operator matrices resulting are unsymmetric (due to the lack of symmetry of the spectral differentiation matrices), non-normal, and not positive definite.  
Lastly due to the ill-posedness of the Poisson-Neumann problem, the operator matrix is rank-deficient, so the system must be made consistent prior to obtaining a solution.  It is worth noting that although the SMPM matrices have several poor properties from a linear algebra perspective, the SMPM is chosen as the discretization here because of its demonstrated prior effectiveness for environmental flow problems, and also because its lack of a weak form makes its implementation relatively straightforward.  

Now we define the SMPM element matrices and inter-element continuity conditions.  In the subsequent discussion all operators (e.g.~$\nabla^2$), element and domain boundaries (e.g.~$\partial\Omega$) and domains are to be understood in the discrete sense; index notation is deprecated in favor of the continuous notation for convenience.  For example, the domain $\Omega$ and the set of discrete points that is used to discretize $\Omega$ are used interchangeably when the intention is clear.
   All derivatives are computed by way of spectral differentiation matrices constructed from the Lagrange polynomial basis.
   Let $Lu = f$ represent the discrete Poisson-Neumann system on $\Omega \subset \R^2$ a domain discretized into an $m_x \times m_z$ element mesh with each element $V_i$ smoothly and invertibly mapped from the master element $[-1,1]\times[-1,1]$.
   On each element a two-dimensional Gauss-Lobatto-Legendre (GLL) grid with $n$ points in each direction is constructed and is used to evaluate the Lagrange interpolant basis and its derivatives by way of spectral differentiation matrices\cite{Costa2000}.  Thus each element contains $n^2$ grid points.
   If $V_i$ and $V_j$ share the $n$ GLL points along one of their  four boundaries, then each element owns a copy of those $n$ GLL nodes in order maintain the discontinuous nature of this method.
   Thus as matrix, $L \in \R^{r\times r}$ is of dimension $r = n^2 m_xm_z$, where $r$ denotes the total number of nodes in the grid $\Omega$.

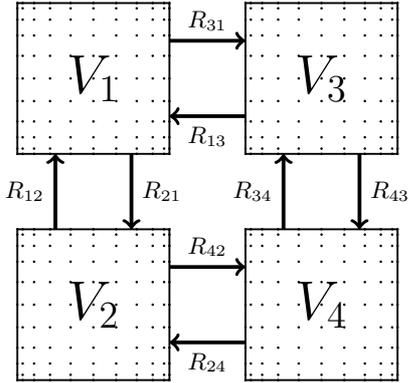
\begin{figure}
	\begin{tikzpicture}
   \label{fig:2x2}

		\draw[line width=0.25mm](0,0) rectangle(2,2);
		\draw[line width=0.25mm](3,3) rectangle(5,5);
		\draw[line width=0.25mm](0,3) rectangle(2,5);
		\draw[line width=0.25mm](3,0) rectangle(5,2);

		\foreach \x in { -1.0, -0.93, -0.78, -0.57, -0.30, 0, 0.30, 0.57, 0.78, 0.93, 1.0 }
			\foreach \y in { -1.0, -0.93, -0.78, -0.57, -0.30, 0, 0.30, 0.57, 0.78, 0.93, 1.0 }
				\draw node[fill,circle,scale=0.1](\x,\y) at (\x+1,\y+1) {};

		\foreach \x in { -1.0, -0.93, -0.78, -0.57, -0.30, 0, 0.30, 0.57, 0.78, 0.93, 1.0 }
			\foreach \y in { -1.0, -0.93, -0.78, -0.57, -0.30, 0, 0.30, 0.57, 0.78, 0.93, 1.0 }
				\draw node[fill,circle,scale=0.1](\x,\y) at (\x+4,\y+1) {};

		\foreach \x in { -1.0, -0.93, -0.78, -0.57, -0.30, 0, 0.30, 0.57, 0.78, 0.93, 1.0 }
			\foreach \y in { -1.0, -0.93, -0.78, -0.57, -0.30, 0, 0.30, 0.57, 0.78, 0.93, 1.0 }
				\draw node[fill,circle,scale=0.1](\x,\y) at (\x+4,\y+4) {};

		\foreach \x in { -1.0, -0.93, -0.78, -0.57, -0.30, 0, 0.30, 0.57, 0.78, 0.93, 1.0 }
			\foreach \y in { -1.0, -0.93, -0.78, -0.57, -0.30, 0, 0.30, 0.57, 0.78, 0.93, 1.0 }
				\draw node[fill,circle,scale=0.1](\x,\y) at (\x+1,\y+4) {};


		\node at (1,1) {\huge{$V_2$}};
		\node at (4,4) {\huge{$V_3$}};
		\node at (1,4) {\huge{$V_1$}};
		\node at (4,1) {\huge{$V_4$}};

		\draw [line width=0.5mm, ->] (2,1.5) -- node[above]{\small{$R_{42}$}}  (3,1.5);
		\draw [line width=0.5mm, ->] (3,0.5) -- node[below]{\small{$R_{24}$}} (2,0.5);

		\draw [line width=0.5mm, ->] (2,4.5) -- node[above]{\small{$R_{31}$}}  (3,4.5);
		\draw [line width=0.5mm, ->] (3,3.5) -- node[below]{\small{$R_{13}$}} (2,3.5);

		\draw [line width=0.5mm, ->] (0.5,2) -- node[left]{\small{$R_{12}$}}  (0.5,3);
		\draw [line width=0.5mm, ->] (1.5,3) -- node[right]{\small{$R_{21}$}} (1.5,2);

		\draw [line width=0.5mm, ->] (3.5,2) -- node[left]{\small{$R_{34}$}}  (3.5,3);
		\draw [line width=0.5mm, ->] (4.5,3) -- node[right]{\small{$R_{43}$}} (4.5,2);

	\end{tikzpicture}
	\caption{A depiction of the logical arrangement of a  $2\times 2$ element spectral multi-domain penalty method (SMPM) grid with $10 \times 10$ Gauss-Lobatto-Legendre points in each element denoted $V_j$ for $j = 1,2,3,4$. The inter-element continuity fluxes are represented with $R_{ij}$ with $i\neq j$. }
\end{figure}

   In the SMPM the weak inter-element continuity condition is of Robin type, and is denoted by the flux $R_{ij} : \partial V_j \longrightarrow \partial V_i$ from element $V_j$ into $V_i$ for $V_i,V_j$ with an adjacent boundary $\partial V_j \cap \partial V_i$ consisting of $n$ grid points.  $R_{ij}$ is defined as
   \begin{align}
      R_{ij} = I + \hat{n}_i \cdot \nabla
   \end{align}
   where $\hat{n}_i: \partial V_i \longrightarrow \R^2$ is the outward pointing normal vector of $\partial V_i$ and $I$ is the identity operator.  A depiction of a $2 \times 2$ element grid with the inter-element fluxes is shown in Fig. 1, in which the elements $V_1,V_2,V_3,V_4$ have been separated to emphasize the discontinuous nature of the SMPM.
Within each element is depicted a $10 \times 10$ GLL grid, and the inter-element fluxes map from the boundary of an element into the boundary of the element its arrow points towards.

The physical boundary conditions are Neumann, and are given on $\partial V_i \cap \partial \Omega$ as $n_i \cdot \nabla$ where $n_i$ is again the outward pointing normal vector.
   Given a function $u$, on an element $V_i$ the residual in the spectral multi domain penalty method is given by the sum of the Laplacian, the inter-element continuity mismatch, and the boundary condition mismatch as
	\begin{align}
      \label{smpm_poisson}
      L_i u_i = \nabla^2 u_i + \tau_i\left( R_{ii} u_i  - \sum_{j\in N(i)} \restr{R_{ij} u_j}{\partial V_i \cap \partial V_j}  \right) + \tau_i \restr{\hat{n}_i \cdot \nabla u_i}{\partial V_i \cap \partial \Omega}  = f_i + \tau_i g_i.
	\end{align}
   Here, $g_i$ is the boundary value of the Neumann boundary condition restricted to element $V_i$, and $N(i)$ is the index set of elements adjacent to $V_i$.  The inter-element continuity, external boundary conditions, and the PDE are all satisfied weakly, since the residual is the sum of these three components.  The penalty parameter $\tau_i$ represents the degree to which the inter-element continuity and boundary conditions are weighted in the residual relative to the PDE, and the optimal choice of $\tau_i$ is determined by stability criteria for hyperbolic problems\cite{Hesthaven1997,Hesthaven1998}, and a heuristic for the Poisson problem \cite{Escobar-Vargas2014}.

\subsection{The Schur complement problem}
\label{schur_assembly}

The domain $\Omega$ is decomposed into $m_x$ many sub-domains $\Omega_i$, each a vertical strip of $m_z$ elements (see, for example, Fig.~\ref{fig_schur_grid})
.  Along each of the $m_x - 1$ interfaces between the sub-domains are $2nm_z$ GLL nodes ($nm_z$ nodes on either side of the interface).
   Denote as $k = 2nm_z(m_x-1)$ the number of interfacial nodes in the domain decomposition, and this set of $k$ interface nodes as $\Gamma$.
   The discrete Poisson operator $L$ (Eq.~\ref{smpm_poisson}) is decomposed into a local term and an inter-subdomain flux term which is used to construct the Schur problem.
   This operator decomposition comprises three operators which are defined now.

   First, denote as $E:\Gamma \longrightarrow \Omega$ the inclusion map that maps from the interfacial grid $\Gamma$ to $\Omega$.
   As a matrix, $E\in \R^{r \times k}$ and is composed of zeros and ones, $E^T$ is the restriction from the full grid to the interface grid $\Gamma$, and $E^T E = I \in \R^{k \times k}$ the identity matrix (naturally $EE^T$ is not an identity matrix).

Second, define an operator $B : \Omega \longrightarrow \Gamma$ that consists of the inter-subdomain Robin boundary fluxes.
$B$ represents all of the the inter-element fluxes $R_{ij}$ for $V_i$ and $V_j$ in different subdomains.
As a matrix, $B \in \R^{k \times r}$, since it computes $I + \hat{n} \cdot \nabla$ within a subdomain using spectral differentiation matrices and assigns it to the interface of its neighbor.

Finally, define the operator $A : \Omega \longrightarrow \Omega$, which represents the part of $L$ that is entirely local to one subdomain.  $A$ consists of the Laplacian part of $L$, the boundary condition mismatch, all inter-element flux terms that are between elements contained in one subdomain, and the $R_{ii}$ terms in Eq.~(\ref{smpm_poisson}).
Since $A$ is entirely local to each subdomain, as a matrix $A\in \R^{r\times r}$ is block-diagonal.
Lastly, note that owing to the Robin-type boundary conditions $R_{ii}$ that are a part of $A$, $A$ represents $m_x$ decoupled homogenous Poisson-Robin boundary value problems, and as such $A$ is invertible and its sparse factorization is easily parallelizable.

These three operators are defined so that their combination yields the SMPM Poisson-Neumann operator
\begin{align}
   Lu = Au + EBu = f.
\end{align}
Using this decomposition of $L$, a Schur problem for the inter-subdomain boundary fluxes $x_S = Bu$ is formed by dividing by $A$ and multiplying by $B$:
	\begin{align}
		\label{operator_decomp}
		Au + EBu = f & \Rightarrow   u + A^{-1}EBu = A^{-1} f \nonumber \\
		& \Rightarrow (I + BA^{-1}E)Bu = BA^{-1}f \nonumber \\
		& \Rightarrow Sx_S = b_S,
	\end{align}
	where the Schur complement matrix is $S = I + BA^{-1}E$, its right hand side $b_S = BA^{-1}f$, and its solution $x_S = Bu$.  Once the boundary fluxes $x_S = Bu$ are obtained by solving the Schur complement system the solution in each subdomain $\Omega_i$ is recovered as
	\begin{align}
		u_i = A_i^{-1}( f_i - (Ex_S)_i).
	\end{align}
	Being block-diagonal, all of the divisions with $A$ are easily parallelized, and so the expensive part of the above is obtaining the solution of the Schur complement system $x_S$.
   In the rest of this paper, the focus is on efficiently obtaining the solution $x_S$ of the Schur complement system.

\subsection{Inconsistency of the Poisson-Neumann system}
\label{inconsistency}
	Prior to obtaining the solution, there remains the important point of dealing with the rank-deficiency of the Poisson-Neumann operator $L$.  The Poisson-Neumann equation is ill-posed in the continuous sense, and so the SMPM operator $L$ is rank-deficient and has non-trivial left and right null spaces of dimension one.  In symmetric discretizations, the kernel vector is the constant vector, but since $L$ is unsymmetric its left and right null spaces are different and only the right null space is constant vector.  To ensure consistency and solvability the right-hand-side vector $f$ is projected out of the left null space of $L$\cite{Pozrikidis2001} and instead of $Lu = f$, the system solved is
	
	\begin{align}
		\label{regularize_poisson}
		Lu =  \tilde{f}
	\end{align}
	where $\tilde{f} = f - u_L u_L^Tf$ is $f$ projected onto the range space of $L$ and $u_L \in \R^r$ is the unique vector with unit norm that satisfies $\norm{ u_L^T L }_2 = 0$.  The solution $u$ then is only known up to an indeterminant additive constant vector.  The rank deficiency of $L$ is inherited by the Schur system, and a relationship between the left null spaces of the Schur system and full Poisson matrix is shown in the following claim.
	
	\begin{cla}  Denote as $u_L \in \R^r$  the left null vector of $L\in \R^{r\times r}$.  Denote as $u_S \in \R^k$ the left null vector of $S \in \R^{k\times k}$.  Then the following relations hold.
		\begin{enumerate}
			\item $u_S = E^T u_L$
			\item $u_L = A^{-T}B^T u_S$
		\end{enumerate}
	\end{cla}
	\begin{proof}
		\emph{See appendix.}
	\end{proof}

	The rank-deficiency of the Schur matrix means another regularization is required to project the Schur right hand side $b_S = BA^{-1}\tilde{f}$ out of the left null space of the Schur system.  Thus the Schur complement system is modified to read
	\begin{align}
	\label{regularize_schur}
		Sx_S = b_S - u_S u_S^T b_S.
	\end{align}
	To summarize, the method for obtaining the solution $u$ to $Lu = f$ is shown in Algorithm \ref{alg:schur}.  The statement SOLVE$(S,b_S)$ in Step \ref{schursolve} is meant to represent the solution of the linear system $Sx_S = b_S$.
	
	\renewcommand{\algorithmicrequire}{\textbf{Input:}}
	\renewcommand{\algorithmicensure}{\textbf{Output:}}
	\begin{algorithm}[H]
		\begin{algorithmic}[1]	
			\REQUIRE $f, u_L, u_S$
			\ENSURE $u$
			\STATE $f \longleftarrow f - u_L u_L^T f$
			\STATE $b_S := BA^{-1}f$
			\STATE $b_S \longleftarrow b_S - u_S(u_S^T b_S)$  \label{reg2}
			\STATE $x_S :=  $ SOLVE($S,b_S$) \label{schursolve}
			\STATE $u \longleftarrow A^{-1}(f - E x_S)$ \label{recover}
		\end{algorithmic}
		\caption{Schur complement method with null space projections.}
		\label{alg:schur}
	\end{algorithm}	

	The second regularization, Step \ref{reg2}, is done to ensure that the Schur system is consistent and thus solvable in Step \ref{schursolve}, but it is not obvious that modifying the right hand side of the Schur system is acceptable from the perspective of solving the Poisson-Neumann equation.  
   How do we know that modifying the right hand side of the Schur system does not give us the wrong answer $x_S$ for Step \ref{recover}?  
   Fortunately, it can be shown that so long as the original Poisson-Neumann system was made consistent (e.g. Eq.~(\ref{regularize_poisson})) the regularization of the Schur problem does not affect the accuracy of the solution of the Poisson equation.
   This is shown in the following claim.
	\begin{cla} Given $L u = \tilde{f} \in \R^r$, if $u_L^T \tilde{f} = 0$, then the error in the solution $u = A^{-1}( \tilde{f} - E x_S )$ recovered from the solution of the regularized Schur complement system $S x_S = b_S - u_Su_S^Tb_S$ is bounded by the error of Schur complement solution: 
		\begin{align}
			\norm{ Lu - \tilde{f}}_2 \leq \norm{Sx_S - (I -  u_S u_S^T)b_S }_2.
		\end{align}
	\end{cla}
	\begin{proof}
		\emph{See appendix.}
	\end{proof}

	Finally, a means to obtain $u_L$ and $u_S$ is required.  Note that first $\norm{u_S^TS}_2 = 0$ means that $u_S$ is a left eigenvector of $S$.  Thus, shifted inverse iteration is used on $S^T - \sigma I$ to find the left kernel vector $u_S$, and then this vector is used to reconstruct the left kernel vector $u_L$ by way of the second relation in Claim 1.  With a small shift of $\sigma$ and a zero eigenvalue $\lambda_n = 0$, the convergence of shifted inverse iteration is of order $( \sigma - \lambda_n )^{-1} = \sigma^{-1}$.  Since $\sigma$ is usually chosen to be close to zero, this iteration converges quickly, often in one or two steps.

\section{Solving the Schur system}
\label{section_solve}
Many preconditioning techniques for the Schur complement system have been proposed\cite{Yamazaki2010a,Cros2002,Carvalho2001,Yamazaki2010}, with an aim towards algorithmic scalability.
By algorithmic scalability is meant the property that an iterative solution method converges independently of the number of degrees of freedom of the underlying grid; in this case it means that the number of GMRES iterations is kept bounded as $m_x$ grows.
For elliptic problems, most preconditioning techniques rely on two preconditioners, a local preconditioner that can be applied in parallel and a coarse global preconditioner to speed across-grid communication of components of the residual.
An example is the two-level additive Schwarz preconditioner in which overlapping block-diagonal components are solved in parallel, augmented with a coarse grid correction to communicate information across the grid (e.g. \cite{Fischer1997,Escobar-Vargas2014}).   In this paper, a non-overlapping block-diagonal/block-Jacobi preconditioner is used, augmented with deflation, to achieve GMRES convergence rates independent of $m_x$.  As $m_x \ll m_z$ in long domain problems, achieving convergence independent of $m_x$ is of crucial importance.

\subsection{Krylov solver implementation}

The iterative solution method used in Step 4 of Algorithm \ref{alg:schur} and throughout this work is the Generalized Minimum Residual Method (GMRES) \cite{Saad2003}.
GMRES was chosen because of its generality; it is capable of solving unsymmetric linear systems like the ones in SMPM, and requires only the ability to compute matrix-vector products of the operator matrix $S$.
Further, even though the Schur matrix inherits singularity from the Poisson matrix, when made consistent the Schur system converges properly in GMRES.
This is first because the (left) null and (right) range spaces of $S$ are orthogonal (i.~e. $N(S) \cap R(S) = \{\phi\}$), a necessary condition for convergence of GMRES \cite{Brown1997} on singular matrices, and second because every principle sub-matrix of $S$ is invertible.
Thus, after projection out of the null space, rank deficiency only plays a role in the final Krylov space which is of course never reached in practice.
However, the non-normality of $S$ ( $SS^T \neq S^T S$) does mean that the conditioning and more generally the spectrum alone is not predictive of GMRES convergence \cite{Greenbaum1996,Meurant2014}.
For the remainder of the discussion, then, iteration counts will be used as the metric for evaluating the performance of the various preconditioning techniques.

In this implementation, the blocks of the Schur matrix, denoted $S_i$ for an interface $\Gamma_i$ between subdomains $\Omega_i$ and $\Omega_{i+1}$, are assembled, stored, and multiplied against vectors in a parallel sparse distributed fashion.
Each such block is of dimension $S_i \in \R^{4nm_z \times 2nm_z}$.
Assembling the Schur matrix takes advantage of the fact that the interfaces are sparse in the full grid for a high-order method and thus assembly and computation of matrix-vector products $\sum_i S_ix_i$ requires less work than the explicit computation of $(I + BA^{-1}E)x$.
So long as storage of the all the blocks $\{S_i\}_{i=1}^{m_x - 1}$ is possible, and the Poisson-Neumann problem is to be solved many times, it is preferable to assemble $S$.

	Finally, note that the implementation of GMRES used here uses Householder reflections for orthogonalization in the Arnoldi process, which requires about three times more operations than the classic modified Gram-Schmidt (MGS) method.  This choice is made because of a loss of orthogonality in the Krylov basis $V_k$ due to round-off error in MGS that scales with condition number \cite{Bjorck1967},
	\begin{align}
		V_k^TV_k = I + \mathcal{O}( \epsilon \kappa_2(S) ),
	\end{align}
	where $\epsilon$ is machine precision and $\kappa_2(S)$ is the condition number of $S$. This is clearly problematic here with a singular $S$.  Walker \cite{Walker1988} shows that by contrast, orthogonalization with Householder transformations yields an orthogonality error of
	\begin{align}
		V_k^TV_k = I +\mathcal{O}( \epsilon)
	\end{align}
	which is a benefit that is worth the extra floating point operations.

\subsection{Non-overlapping block-Jacobi preconditioner}
\label{sec:bj}

\begin{figure}[h]
		\begin{center}
		\begin{picture}(200,100)
			\put(-100,0){ \includegraphics[width=0.875\textwidth]{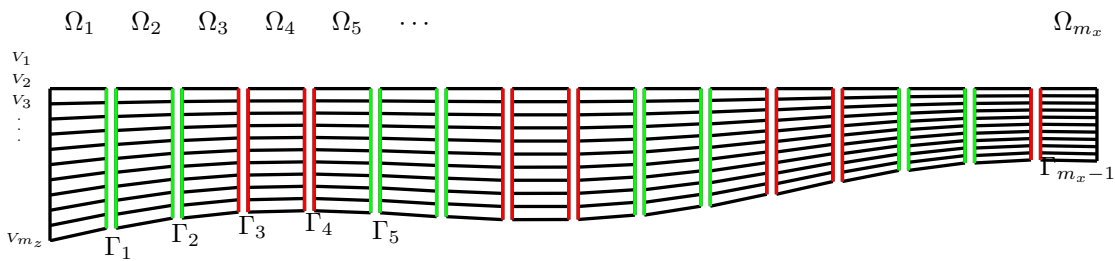}}
			\put(-90,80){$\Omega_1$}
			\put(-65,80){$\Omega_2$}
			\put(-40,80){$\Omega_3$}
			\put(-15,80){$\Omega_4$}	
			\put(10,80){$\Omega_5$}				
			\put(35,80){$\cdots$}	
			\put(280 , 80){$\Omega_{m_x}$}		
			\put(-75,-5){$\Gamma_1$}	
			\put(-50,0){$\Gamma_2$}			
			\put(-25, 3){$\Gamma_3$}			
			\put(-0,4){$\Gamma_4$}	
			\put(25,1){$\Gamma_5$}		
			\put(275,25){$\Gamma_{m_x - 1}$}																														
			\put(-110,68){\tiny{$V_1$}}
			\put(-110,60){\tiny{$V_2$}}
			\put(-110,52){\tiny{$V_3$}}
			\put(-108,38){\tiny{$\vdots$}}			
			\put(-112,0){\tiny{$V_{m_z}$}}
		\end{picture}
		\end{center}
		\caption{ A sample domain with $m_x = 17$ and $m_z = 10$ elements in $x$ and $z$ respectively.  Each vertical strip of $10$ elements is a sub-domain $\Omega_j$, and the collection of interfaces $\{\Gamma_k\}_{k=1}^{m_x -1}$ represents the grid for the Schur complement problem. The red/green interfaces correspond to the red/green blocks in the sparsity pattern of the block-Jacobi preconditioner shown in Fig.~\ref{fig_schur}.}
		\label{fig_schur_grid}
	\end{figure}

For sparse matrices with non-zeros clustered around the diagonal, computing the inverse of  blocks along the diagonal separately can be a useful preconditioning technique.
So-called block-Jacobi preconditioners have been shown to be effective for the Schur complement of elliptic operators \cite{Couzy1995}, especially when combined with coarsened-grid preconditioners \cite{Manna2004, Pavarino2000, Pasquetti2006}.
In the context of solving the Poisson equation, $A^{-1}$ as defined in Section \ref{schur_assembly} represents a block-Jacobi preconditioner since $A$ represents decoupled local problems, each a Poisson-Robin boundary value problem.
For preconditioning the Schur matrix $S$, a block-Jacobi preconditioner can also be assembled in which a single block represents the coupling between interfaces $\Gamma_i$ and $\Gamma_{i+1}$ in the Schur grid $\Gamma$.  Each such block represents the coupling between four interfaces: the two bounding a subdomain $\Omega_i$ on the left and the right, and their adjacent interfaces in $\Omega_{i-1}$ and $\Omega_{i+1}$.  This is pictorially depicted in Fig.~\ref{fig_schur_grid}, in which every four consecutive interfaces belonging to one block in the block-diagonal preconditioner are shown in the same color (either red or green), and represent the nodes belonging to $\Gamma_{i}$ and $\Gamma_{i+1}$.   

\begin{figure}[h]
		\begin{center}
		\includegraphics[width=0.50\textwidth]{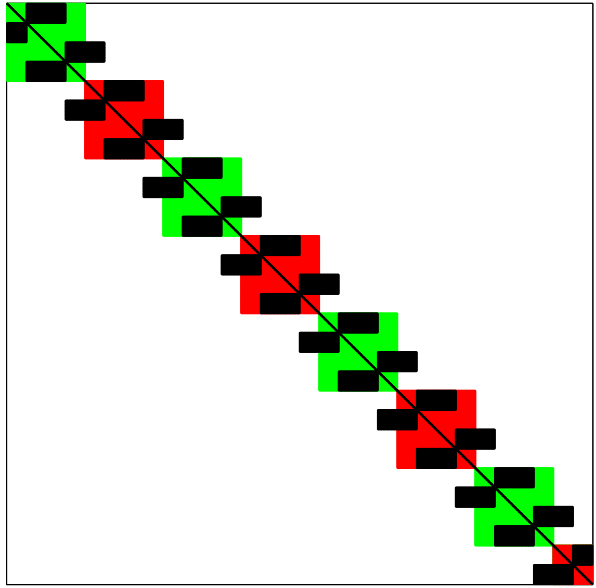}
		\end{center}
		\caption{ The Schur matrix shown in black corresponding to a grid with $n = 5$, $m_x = 17$, and $m_z = 10$ shown in Fig.~\ref{fig_schur_grid}, with the blocks $M_i$ of the non-overlapping block-Jacobi preconditioner $M$ in red and green.  Each block-Jacobi block represents the coupling between four of the interfaces shown in Fig.~\ref{fig_schur_grid}, in which the red/green interfaces correspond to the red/green blocks shown here.}
		\label{fig_schur}
	\end{figure}

To graphically show the relationship between the grid and the blocks in the block-Jacobi preconditioner, the sparsity pattern of the Schur matrix corresponding to the grid in Fig.~\ref{fig_schur_grid} is depicted in black in Fig.~\ref{fig_schur}.   The dimensions of the grid are $n = 5, m_z = 10, m_x = 17$.
Overlaid on the sparsity pattern of $S$ in Fig.~\ref{fig_schur} are the blocks of the block-Jacobi preconditioner used here, shown in alternating colors of green and red that correspond respectively to the green and red interfaces they represent in Fig.~\ref{fig_schur_grid}.
Since each interface consists of $nm_z$ grid points, the dimension of each square block is $4nm_z$, except for possibly the last such block which may be smaller if $m_x$ is odd.
Denoting the block-Jacobi preconditioner matrix as $M$, the preconditioned Schur system that is solved with GMRES is
\begin{align}
   SM^{-1}x'_S = b_S
\end{align}
and the solution is obtained by a final division by $M$
\begin{align}
   x_S = M^{-1}x'_S.
\end{align}
Since $M$ is explicitly block-diagonal (i.e. any non-zeros of $S$ coupling the blocks of $M$ are ignored in the factorization of $M$), divisions by $M$ can be computed efficiently in parallel.  


	\subsection{Deflation}
   \label{sec:def}
Working in tandem with other preconditioners, deflation methods aim to accelerate the convergence of Krylov methods by eliminating (or ``deflating'') components of the residual within a chosen subspace.  The subspace is usually chosen to be a span of approximate eigenvectors of the operator corresponding to slowly converging eigenvalues.  Thus, the problematic eigenvalues are solved directly using a coarsened version of the operator, and the remaining components of the residual are eliminated by a Krylov solver.  The first deflation methods were used to accelerate the convergence of the conjugate gradient method for elliptic boundary value problems (e.g. \cite{Nicolaides1987}), and thus were limited to symmetric positive definite matrices.  The extension to unsymmetric matrices \cite{Erlangga2008} and in particular applications related to domain decomposition with preconditioning \cite{Vuik2006} and fluid flow \cite{Tang2012} make deflation a good candidate for augmenting block-Jacobi preconditioning in our Schur problem.  Comparisons of coarse grid correction with deflation show that deflation methods combined with preconditioning yield lower condition numbers irrespective of the choice of the coarsening operators for symmetric positive-definite systems\cite{Nabben2004}, and deflation methods have been elsewhere experimentally favorably compared in the context of two-level preconditioners with multi-grid and domain decomposition approaches \cite{Tang2009,Tang2010}.  A more exhaustive summary and analysis of deflation methods can be found in Ref. \cite{Gaul2013}.

	Here, to augment the block-Jacobi preconditioner described in Section \ref{sec:bj}, a deflation method is used as the coarse-grid correction method, following the procedure in Ref.~\cite{Erlangga2008}.
  The deflation vectors are chosen to be a set of $d$ column vectors $Z \in \R^{k\times d}$ where $d \ll k$, and $k = \textrm{dim}(S)$.
   These deflation vectors are chosen to be discrete indicator vectors on each of the interfaces $\{\Gamma_j\}_{j=1}^{m_x - 1}$ that form the Schur grid $\Gamma$.
   The $i$-th entry in the $j$-th such vector is given by
	\begin{align}
		(z_j)_i = \left\{
		     \begin{array}{ll}
		       1 & :  \textrm{ if $x_i \in \Gamma_j$ } \\
		       0 & : \textrm{ if $x_i \notin \Gamma_j$ }
		     \end{array}
		     \right\},
	\end{align}
	thus each vector is active on one interface in the Schur grid.  The matrix of these vectors  $Z = [ z_1, z_2, \cdots, z_d ]$  defines a coarse version of the Schur problem, $C = Z^TSZ \in \R^{d\times d}$, and two projections
		\begin{align}
			P &= I - S ZC^{-1}Z^T \label{projection1}	 \\
			Q &= I - ZC^{-1}Z^TS
			\label{projection2}
		\end{align}
		each of size $\R^{k\times k}$.
      The operator $Z$ defines a one-dimensional coarse grid of dimension $d$ that is composed of the mean $x$-coordinates of each of the $m_x-1$ interfaces $\Gamma_j$.
      As a matrix $Z^T \in \R^{ d \times k}$ is a contraction operator that maps grid functions on the Schur grid to the coarse grid.
      The mapping $y = Z^Tx$ sums up the values of $x$ along each interface $\Gamma_j$ and stores them in $y_j$.
      The mapping $x = Zy$ is the prolongation operator which assigns to all entries of $x$ on $\Gamma_j$ the value $y_j$.
      Note that $Z^TZ = nm_z I \in \R^{m_x - 1 }$.
      $ZZ^Tx$ averages the values of $x$ over each interface $\Gamma_j$.   The intuition behind the projections $P$ and $Q$ is that they project out of the subspace on which $ZC^{-1}Z^T$ is a good approximation of the left (in the case of $Q$) or right (in the case of $P$) inverse of $S$.  Thus the projections map onto the complement of the subspace on which the coarse matrix $C$ approximates the Schur matrix $S$ well.

		Departing from the method in Ref. \cite{Erlangga2008}, the division by $C$ is regularized to address $C$'s  rank-deficiency inherited from $S$.
      The projections $P$ and $Q$ are also regularized in the same way.
      Note that the left null space of $C$ is one dimensional, and is spanned by $Z^Tu_S$.
      Denote this vector $u_C = Z^Tu_S$ and write the regularized projections as
		\begin{align}
			P &= I - S Z (C \backslash (Z^T - u_Cu_C^TZ^T)) \\
			Q &= I -  Z (C \backslash (Z^T - u_Cu_C^TZ^T))S.
		\end{align}
		Deflation proceeds by noting that the solution of the Schur complement problem $Sx_S = b_S$ can be decomposed into
		\begin{align}
			x_S = (I - Q)x_S + Qx_S.
		\end{align}
		Then, the first term is just $ZC\backslash( Z^T - u_C u_C^TZ^T)(b_S - u_Su_S^Tb_S)$ which can be computed directly since $C$ is small, and the second term is obtained by performing GMRES on the deflated and right-preconditioned system $PSM^{-1}x_S = P(b_S - u_Su_S^Tb_S)$, and then post-multiplying by $Q$, finally assembling the solution as
		\begin{align}
         \label{deflated_schur}
			x_S &= ZC\backslash( Z^T - u_C u_C^TZ^T)(b_S - u_Su_S^Tb_S)  \nonumber \\ &+ QM^{-1}\textrm{GMRES}(PSM^{-1},P(b_S - u_Su_S^Tb_S)).
		\end{align}
		Because $P$ projects out of the coarse space, the GMRES solution of $PSM^{-1}x_S = P( b_S - u_Su_S^Tb_S)$ minimizes only the component of the residual that cannot be well-approximated by the coarse solution.
      Again note that although the deflated operator $PSM^{-1}$ is rank-deficient by virtue of the rank-reducing projection $P$, because the the right-hand-side vector is similarly projected into the column space the GMRES algorithm converges properly.
      This formulation of deflation-augmented right-preconditioning is an extension of the work in Ref.~\cite{Erlangga2008} to a rank-deficient matrix.

      With reference to implementation, the matrices $Z^T$ and $Z$ are very easily parallelizable, and since $P$ and $Q$ are never explicitly formed, the storage overhead for deflation is the storage of the coarse matrix $C$ which is a one-dimensional tridiagonal finite-difference Poisson matrix of dimension $m_x-1$ whose solution is obtained by the Thomas algorithm.
      The applications of $P$, $S$, and $M^{-1}$ required in each iteration of GMRES are conducted in a distributed, sparse, and MPI-parallel fashion, with the exception of the coarse solution, which requires global communication between all MPI ranks.

	Finally, note that the two-level additive Schwarz method can also be written using the deflation vectors $Z$ defined above as coarsening operations. Denoting the coarse problem again as $C = Z^TSZ$, write the two-level additive Schwarz	preconditioned system as 
	\begin{align}
		S(M^{-1} + ZC\backslash( Z^T - u_C u_C^T Z^T))x_S = b_S - u_Su_S^{T}b_S.
	\end{align}
	As defined above, the two-level additive Schwarz preconditioner will be used to compare against the performance of the deflation method.
   Notice that although they appear similar, the additive Schwarz method only requires one application of the Schur matrix $S$ per Krylov iteration whereas the deflation method requires two (a second embedded in the definition of $P$).  Furthermore, unlike the deflation method does with its projection $P$, additive-Schwarz preconditioning doesn't explicitly project out of the coarse space.  Thus while the Krylov solver is aided by the coarse matrix $C$, the residual it minimizes still contains components in the coarse space.
   Finally notice that right-preconditioning has been used throughout.
   This is so that the the tolerance achieved on the preconditioned system bounds the tolerance on the unpreconditioned system, and not the other way around as when left preconditioning is used.

	For completeness, in Algorithm \ref{alg:seq} is the algorithmic summary of the deflation method in which the notation GMRES$(A,b)$ is intended to represent the solution of a linear system $Ax = b$ with GMRES.  The algorithm expects to be given the vectors spanning the left null spaces of $C,S,$ and $L$ ($u_C$,$u_S$,$u_L$ respectively) and returns the solution $x$ that makes $Lx - b$ small.  Naturally there are input parameters for GMRES (desired tolerance, maximum iterations, etc.) that are left out in the description below.
	\renewcommand{\algorithmicrequire}{\textbf{Input:}}
	\renewcommand{\algorithmicensure}{\textbf{Output:}}
	\begin{algorithm}[H]
		\begin{algorithmic}[1]	
			\REQUIRE $b, u_L, u_S, u_C$
			\ENSURE $x$
			\STATE $b \longleftarrow b - u_L u_L^T b$
			\STATE $b_S := BA^{-1}b$
			\STATE $b_S \longleftarrow b_S - u_S(u_S^T b_S)$ 
			\STATE $x_1 :=  $ GMRES($PSM^{-1},Pb_S$)
			\STATE $x_1 \longleftarrow QM^{-1} x_1$
			\STATE $x_2 := Z^T b_S - u_Cu_C^TZ^Tb_S$ 
			\STATE $x_2 \longleftarrow ZC\backslash x_2$
			\STATE $x := x_1 + x_2$
			\STATE $x \longleftarrow A^{-1}(b - E x)$
		\end{algorithmic}
		\caption{Block-Jacobi Preconditioned and Deflated Schur Method}
		\label{alg:seq}
	\end{algorithm}

\section{Performance}
\label{sec:2d}

In this section a performance comparison of four different preconditioning methods for solving the Schur complement problem using GMRES is made.
The four methods are 
\begin{enumerate}
   \item no preconditioning,
   \item block-Jacobi preconditioning,
   \item deflation-augmented block-Jacobi preconditioning,
   \item two-level additive-Schwarz preconditioning.
\end{enumerate}
These are all described in Section \ref{section_solve}, and are summarized in Table \ref{table_methods}.
The basis for comparison will be the number of iterations for achieving a GMRES relative tolerance of $10^{-10}$, as well as the wall-clock time taken for achieving that tolerance.
Both iteration count and computation time are relevant since each of the four methods has a different per-Krylov iteration computational cost.
In all cases, the right hand side used is a randomly generated vector drawn from a uniform distribution on $[0,1]$, and the measurement of iteration count and computation time is averaged over the solution of ten such right hand sides.
The initial guess is always the vector of zeros, and the solution $x_S$ is always checked to verify that
\begin{align}
   \frac{ \norm{Sx_S - b_S}_2}{\norm{b_S}_2} < 10^{-10}.
\end{align}

Two classes of problems are studied to show GMRES convergence independent of grid properties.
First, to show that the block-Jacobi preconditioner alone eliminates the dependence on domain aspect ratio, the Poisson problem is solved on a series of increasingly leptic (stretched) grids.
Second, to show that the projections in the deflation method eliminate the dependence on $m_x$, a series of grids with increasingly many $x$ elements is constructed on which the Poisson equation is solved.
In the second set of problems, deflation is compared with two-level additive Schwarz preconditioning to demonstrate that deflation requires both fewer GMRES iterations and less wall-clock time.
All computations are performed in an MPI-parallel Fortran code and are executed on a 64-core AMD Opteron computer using the AMD Core Math Library for all basic linear algebra tasks.

	\begin{table}[h]
	\caption{A description of the four methods compared in this section.  GMRES$(A,b)$ represents a GMRES solution of the linear system $Ax = b$.}
	\begin{center}
	\begin{tabular}{lll}
		\bf{Name} & \bf{Abbreviation} & \bf{Description} \\
		No Preconditioning &	Schur & GMRES($S, \tilde{b}_S$) \\
		Block-Jacobi  &	BJ-Schur & GMRES( $SM^{-1},\tilde{b}_S$) \\
		Deflation and Block-Jacobi  &	DBJ-Schur & GMRES($PSM^{-1}, P\tilde{b}_S$) \\
		Two-Level Additive Schwarz  &	2LAS  & GMRES($S(M^{-1} + ZC\backslash( Z^T - u_C u_C^T Z^T)), \tilde{b}_S$ )
	\end{tabular}
	\label{table_methods}
	\end{center}
	\end{table}


\subsection{Increasing domain aspect ratio}

The goal of this section is to demonstrate that use of the block-Jacobi preconditioner on the Schur matrix eliminates the dependence of GMRES convergence on the domain aspect ratio.
To show this, the Poisson-Neumann problem was solved on a series of domains with constant height of $l_z = 10$ m and incrementally increasing width from $l_x = 10$ m to $l_x = 5000$ m.
These domains are discretized with $10$ elements in both the horizontal and the vertical, yielding an element aspect ratio $\eta = (l_x/m_x)/(l_z/m_z)$ that grows from $\eta = 1$ to $\eta = 500$.

On this set of grids the Schur problem is assembled and solved with and without block-Jacobi preconditioning.
Shown in Fig.~\ref{fig_residual}(b) are the residual histories of the unpreconditioned GMRES solver for five such grids with aspect ratios increasing from $\eta = 1$ to $\eta = 500$.
Notice that as the aspect ratio increases, there is an increasingly long period of slow decay of the GMRES residual before it quickly decays to the desired tolerance.
As $\eta$ increases to 500, the number of GMRES iterations grows from under about 80 to 1500, a twenty-fold increase that is expected due to the ill-conditioning of the spectral differentiation matrices that worsens as the element aspect ratio grows\cite{Hesthaven2007}.

\begin{figure}[h]
		\begin{center}
		\includegraphics[width=0.95\textwidth]{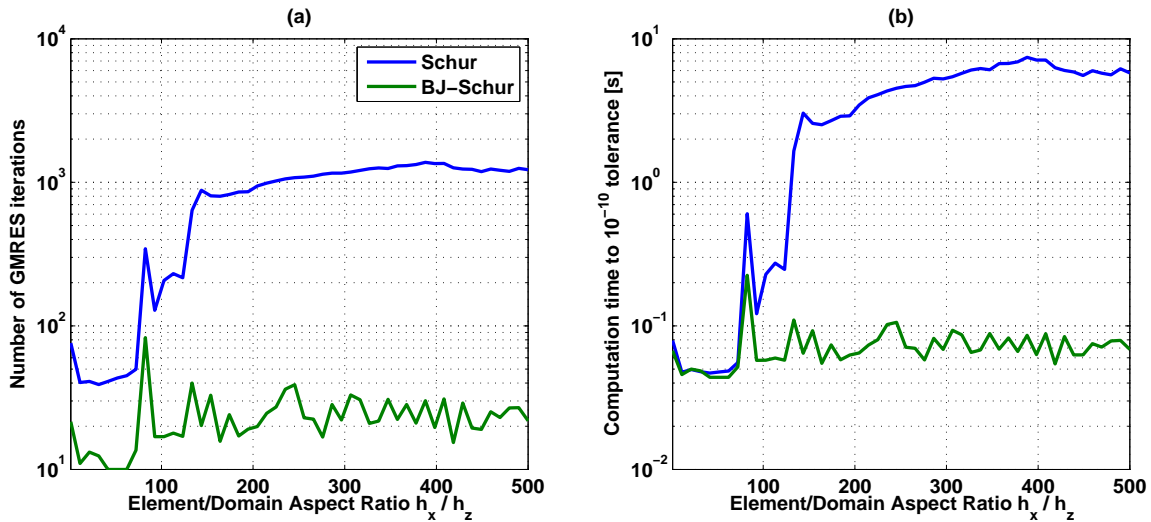}
		\end{center}
		\caption{\emph{Left}: Number of GMRES iterations to solution as a function of element (and domain) aspect ratio for a $10 \times 10$ subdomain grid with $100$ points per element. \emph{Right}: Time to solution the same problem as a function of element (and domain) aspect ratio.  This was benchmarked with an MPI-parallel Fortran code on 5 MPI ranks.}
		\label{fig_aspect_fortran}
	\end{figure}

As shown in Fig.~\ref{fig_aspect_fortran}(a), and in contrast to the unpreconditioned Schur method, the same set of problems solved with block-Jacobi preconditioning shows no dependence on aspect ratio, converging to the tolerance in roughly $20$ GMRES iterations for all values of $\eta$.
Notice though that the iteration count in the unpreconditioned Schur method does not grow indefinitely; it plateaus after $\eta > 150$ at its final value of $\sim 1000$.
This behavior can be explained by understanding that for $\eta \approx 150$ and larger, all horizontal wavenumbers supported by the grid are longer than the supported vertical wavenumbers.
   As such, the discrete eigenvalue spectrum of the Poisson operator is separated into two distinct sets, a condition that defines leptic grids.
   This separation of eigenvalues presents a difficulty for Krylov subspace methods as has been demonstrated previously \cite{Santilli2015}, but once the horizontal and vertical eigenvalues are separated entirely further distortion of the grid only minimally affects the convergence properties of GMRES.
   Finally, both methods show a peak in iterations near $\eta = 75$, a result that is not well understood.

Commensurate with the difference in iteration count, block-Jacobi preconditioned Schur is significantly faster than the unpreconditioned Schur method.
The wall-clock computation time for both methods is shown in Fig.~\ref{fig_aspect_fortran}(b), which demonstrates that with the use of the block-Jacobi preconditioner a solution is obtained in $\mathcal{O}(10^{-1})$ seconds for all aspect ratios. 
By contrast, the time to solution for the unpreconditioned Schur method grows from $\mathcal{O}(10^{-1})$ to $\mathcal{O}(10^2)$ seconds as the aspect ratio increases.
Lastly note that the results for the deflated and two-level additive Schwarz solvers are not shown, since for a small $10 \times 10$ element grid the coarse-grid correction does not improve performance significantly, and can often degrade performance in parallel due to the additional inter-processor communication required in solving the coarse matrix.

\subsection{Increasing the number of $x$ elements}
\begin{figure}[H]
		\begin{center}
		\includegraphics[width=0.95\textwidth]{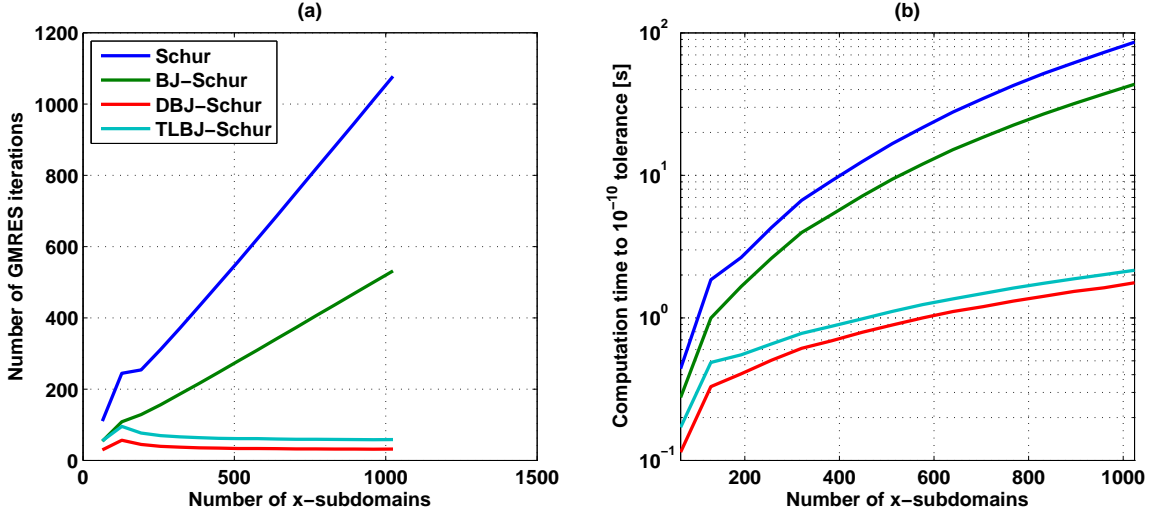}
		\end{center}
		\caption{\emph{Left}: Number of GMRES iterations to solution for a constant aspect ratio as a function of the number of $x$ subdomains. \emph{Right}: Time to solution for a constant aspect ratio as a function of the number of $x$ subdomains. This simulation was benchmarked with $(n, m_z) = (10, 10)$ fixed, and a GMRES tolerance of $10^{-10}$ on 16 processors in an MPI-parallel Fortran code.}
		\label{fig_mx_fortran}
	\end{figure}

   Next is the study of the performance of the various GMRES preconditioning methods as the number of $x$ elements is increased.
   This study gives an insight into how the four methods (unpreconditioned, block-Jacobi preconditioned, two-level additive Schwarz, and deflation) will perform on environmental flow problems in which not only is the domain aspect ratio large, but the number of elements in the long direction is also large ($m_x \gg m_z$).
   In this study, $n = 10$ and $m_z = 10$ were fixed, and the Poisson-Neumann problem is solved for discretizations with $m_x = 64$ to $m_x = 1024$ on domains with increasing length $l_x$.
   The size of these grids grows from $64000$ when $m_x = 64$ to $1024000$ when $m_x = 1024$.
   For each such grid, the Schur problem was assembled and solved with the four preconditioning techniques described at the beginning of this section.
   Each computation was performed ten times and the iteration counts and timings were averaged over these 10 trials.
    Fig.~\ref{fig_mx_fortran}(a) shows the resulting number of iterations required to achieve a tolerance of $10^{-10}$ for each of the four solvers as a function of the number of $x$ subdomains $m_x \in [64, 1024]$.
   Notice that the deflation/block-Jacobi preconditioner achieves GMRES convergence independent of $m_x$ and converges within approximately $30$ iterations.
   Convergence independent of $m_x$ is vital for reliable performance in applications in long domains.
   In all cases, deflation with block-Jacobi preconditioning required the fewest iterations, and converged within approximately $30$ iterations independent of the problem size.
   The two-level additive Schwarz method also shows convergence independent of problem size but takes twice as many GMRES iterations as deflation in all cases.

	Since the cost of $K$ GMRES iterations is $\mathcal{O}(K^3)$, and iterations grow in the unpreconditioned Schur method linearly with $m_x$, it is expected that the any coarse-grid corrected method cubically outperforms the unpreconditioned Schur method in solution time as $m_x$ becomes large.
   In practice, the speedup is closer to quadratic than cubic in $m_x$ because the solution of the coarse problem requires far more communication than the Schur method alone, and there is also some ancillary computation to set up the coarse solve in $C$.

   A comparison of the four methods' time-to-solution is presented in Fig.~\ref{fig_mx_fortran}(b), which shows that the speedup in deflation over block-Jacobi preconditioned Schur is $\approx \mathcal{O}(m_x^2)$.  The two-level additive Schwarz method shows similar quadratic in $m_x$ speedup relative to the Schur method, but is always slower than the deflation method since it takes about twice as many GMRES iterations to converge.  The difference in time between the deflated and additive-Schwarz methods is smaller than the factor of two difference in GMRES iterations because each GMRES iteration of the deflation method requires two applications of the $S$ matrix due to the additional $S$ application in the $P$ projection; the two-level additive Schwarz method requires only one $S$ matrix-vector multiply.  Nevertheless, deflation augmented preconditioning outperforms two-level additive Schwarz preconditioning in terms of wall-clock time for all of the grids studied here by about $25\%$.

\begin{figure}[H]
		\begin{center}
		\includegraphics[width=0.95\textwidth]{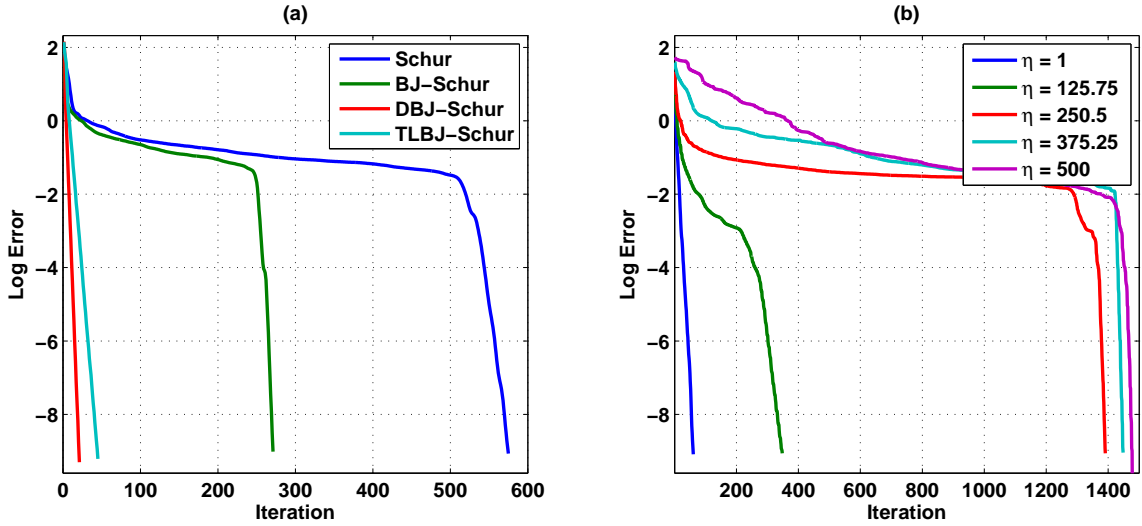}
		\end{center}
		\caption{ \emph{Left}: The GMRES relative error as a function of iteration for a $10 \times 512$ subdomain grid with $n =10$ for each of the four methods compared.  \emph{Right}:  The GMRES relative error as a function of iteration for five different aspect ratios ($\eta = h_x/h_z$) with the unpreconditioned Schur method on a $10 \times 10$ subdomain grid with $n = 10$. }
		\label{fig_residual}
	\end{figure}

		Fig~\ref{fig_residual}(a) shows a representative example of the residual history for the case where $m_x = 512$.  The residual histories of both the block-Jacobi and unpreconditioned cases show an initial period of slow convergence followed by a rapid decay in the residual.  This slow convergence is associated with the first few singular values which represent low-frequency components within the Schur matrix whose elimination requires across-grid communication.  This initial slow convergence behavior is eliminated by either the coarse-grid solve in two-level additive Schwarz or by the deflation vectors in the deflation preconditioner.  However, it is clear that the deflation method yields better GMRES performance since the residual decays at a rate twice that of the two-level additive Schwarz residual.

	To summarize, although both deflation and additive-Schwarz methods demonstrate algorithmic scalability in $m_x$, it is always preferable to use the deflated Schur method over the two-level additive Schwarz method, both from the perspective of minimizing storage (in the number of Krylov basis vectors required) as well as minimizing computation time. This confirms for the Schur complement problem the results established for symmetric positive definite discretizations of the Poisson-Dirichlet operator \cite{Nabben2004}.   The deflation method takes about half as many iterations, and about $25\%$ less time.  The results shown in Fig.~\ref{fig_mx_fortran} are quantified in Table \ref{table_time}, along with data showing the amount of setup time required to assemble the Schur complement matrix, its coarse version, and factor its preconditioner.  Setup takes about 100 times more time than solution, meaning that for the deflation method to be worthwhile, one should require the Poisson-Neumann system to be solved many more than 100 times.  For the vast majority of time-dependent Navier-Stokes simulations, the Poisson-Neumann problem is to be solved many more than $\mathcal{O}(10^2)$ times; its solution is required once per time-step, a total of $\mathcal{O}(10^4-10^5)$ many times \cite{Diamessis2005,Escobar-Vargas2014}.

     \begin{table}[H]
   	\begin{center}
	\caption{Comparison of iterations and computation time to solution in seconds on long domains with randomly generated right hand sides for the four methods averaged over ten trials of randomly generated right hand sides.  In all cases $(n,m_z) = (10, 10)$, the GMRES tolerance was $10^{-10}$, and all were run on 32 processors.}
	\label{table_time}
	\begin{small}
    	\scalebox{0.8}{	
        \begin{tabular}{ r r  r | c  c | c  c | c  c | c  c | }        
            \cline{4-11}
             & & &  \multicolumn{2}{|c|}{Schur} & \multicolumn{2}{|c|}{Block Jacobi }  & \multicolumn{2}{|c|}{Deflation}& \multicolumn{2}{|c|}{TL Schwarz}\\
            $m_x$  & Grid Points & Setup Time & Iter. & Time & Iter. & Time & Iter. & Time & Iter. & Time \\
            \hline 
            \multicolumn{1}{|r|}{64}    & $6.4 \times 10^4$       & 4.24e1 & 110.5    & 4.40e-1 & 54.1    & 2.76e-1  & 29.6   & 1.14e-1  & 54.9   & 1.71e-1  \\
            \multicolumn{1}{|r|}{128}  & $1.28 \times 10^5$     & 8.49e1 & 244.4    &  1.85e0 & 108.6  & 9.97e-1  & 56.7   & 3.29e-1  & 95.4   & 4.87e-1 \\
            \multicolumn{1}{|r|}{256}  & $2.56 \times 10^5$     & 1.68e2 & 311.1    & 4.30e0  & 156.2  & 2.61e0    & 39.3  & 5.04e-1   & 69.7  & 6.56e-1\\
            \multicolumn{1}{|r|}{512}   &$5.12 \times 10^5$     & 3.27e2 & 557.7    & 1.67e1 &  278.6  & 9.40e0   & 33.4   & 8.92e-1   & 61.2  & 1.11e0 \\
            \multicolumn{1}{|r|}{1024} & $1.024 \times 10^6$  & 6.60e2 & 1077.7  & 8.61e1 &  531.8  & 4.35e1   & 31.9   & 1.76e0    & 58.6  & 2.15e0\\
            \hline
        \end{tabular}
        }
        \end{small}
        \end{center}
  \end{table}		

    	\section{Extension to 3D}
      \label{sec:3d}
   All of the preceding discussion has been strictly two-dimensional.  
   A straightforward extension of the deflation method and domain decomposition to the three-dimensional problem is possible if it is assumed that the solution $u(x,y,z)$ is periodic in the $y$ direction, and that the domain is constant in the $y$ direction.
   The first assumption allows for a Fourier discretization in the transverse $y$ direction, and the second assumption allows for the re-use of the domain decomposition outlined in Section \ref{schur_assembly} in the three-dimensional case.
   In this section is discussed the construction and solution of the three-dimensional Poisson-Neumann Schur complement problem with a Fourier transverse discretization.

Denoting as $\Omega \subset \R^2$ a two-dimensional domain, denote as $\Omega' = \Omega \times [0,l_y] \subset \R^3$ the domain on which the three-dimensional Poisson equation is to be solved.
Assume that the boundary conditions in the third dimension $y \in [0,l_y]$ are periodic, implying a periodic solution $u(x,y=0,z) = u(x,y=l_y,z)$, and write the Poisson problem as
	\begin{align}
		\nabla^2 u &= f  \textrm{    on    } \Omega' \nonumber \\
		n\cdot \nabla u &= g \textrm{    on   } \partial\Omega \times [0,l_y] , \nonumber \\
      u(x,0,z) &= u(x,l_y,z) \textrm{   on   } \Omega \backslash \partial \Omega \times \{0,l_y\}.
		\label{poisson3d}
	\end{align}
This type of extension from two to three dimensions facilitates a Fourier discretization in the third ($y$) dimension, and has been used previously \cite{Diamessis2005,Karniadakis1989} in the context of solving the Navier-Stokes equations on problems and domains amenable to transverse periodicity of the solution.
Admittedly requiring periodicity in the third dimension is a limitation and thus this is not a general three-dimensional formulation.
However, as is the case in oceanic and atmospheric modeling of turbulence, only a transect of the large physical domain is discretized computationally, and the transverse direction is assumed to be the direction in which turbulence is statistically homogeneous and thus periodicity can be safely assumed.  Therefore while restrictive, a periodic transverse direction still has significant applicability.

\subsection{Construction of the three-dimensional Schur problems}
   Starting with the Poisson problem (Eq.(\ref{poisson3d})), force periodicity of $u(x,y,z)$ in $y$ by taking the expansion of $u$ in the Fourier basis,
   \begin{align}
	u(x,y,z) = \sum_{j =0}^{m_y/2 - 1} \hat u_j(x,z) e^{ik_j y} ,
   \end{align}
   where $k_j = 2\pi j / h_y$ is the transverse wavenumber, $h_y = l_y/m_y$ spacing of the uniform grid in the transverse direction, $m_y$ the number of grid points in the transverse direction, and $\hat{u}_j(x,z) \in \C$ the Fourier coefficients.
   Substituting this Fourier expansion into Eq.~(\ref{poisson3d}), for each $k_j$ wavenumber a two-dimensional Helmholtz equation in $x$ and $z$ is obtained,
	\begin{align}
		\label{fourier_poisson}
		\nabla^2 \hat u_j(x,z) - k_j^2 \hat u_j(x,z) = \hat f_j(x,z),
	\end{align}
	where $\hat u_j, f_j$ are the $k_j$-th wavenumber components of the Fourier transforms of $u, f$ along the $y$ direction,
   \begin{align}
      \hat u_j &= \hat u( x,k_j,z )  \\
      \hat f_j &= \hat f( x,k_j,z ),
   \end{align}
   and $\hat u(x,k_y,z) = \mathcal{F}_y u(x,y,z)$, $\hat f(x,k_y,z) = \mathcal{F}_y f(x,y,z)$, where $\mathcal{F}_y$ is the discrete Fourier transform in $y$.

   A complete description of the addition of a third Fourier dimension for the spectral multidomain penalty and spectral element methods can be found in Refs.~\cite{Diamessis2005} and ~\cite{Karniadakis1989} respectively, but the important fact is that
	all of the Schur complement methodology described thus far applies directly to each $k_j$ wavenumber in Eq.~(\ref{fourier_poisson}).
   To see this, write the discrete version of Eq.~(\ref{fourier_poisson}) for the $j$-th wavenumber as
	\begin{align}
		( L - k_j^2 I ) \hat u_j = \hat f_j
	\end{align}
  	where now $u_j, f_j \in \C^{k}$ and $L$ is the SMPM Poisson-Neumann operator as in Eq.~(\ref{operator_decomp}).  Using the decomposition of $L = A + EB$ again as in Eq.~(\ref{operator_decomp}) write
	\begin{align}
		(A - k_j^2 I) \hat u_j + EB \hat u_j = \hat f_j.
	\end{align}
   This is a Helmholtz equation with block-diagonal component $A - k_j^2I$ and off-diagonal component $EB$.
   Denoting the shifted block-diagonal matrix $A(k_j) = A - k_j^2I$, the Schur complement problem for wavenumber $k_j$ is given by
	\begin{align}
      \label{transverse_schur}
		S(k_j) = I + BA(k_j)^{-1}E,
	\end{align}
   analogous to the unshifted case considered in Section \ref{schur_assembly}.

	There are two important things to notice about Eq.~(\ref{transverse_schur}).
   First note that while the $B$ matrix can depend on $\tau$, the penalty parameter, and in turn $\tau$ can depend on the shift $k^2_j$\cite{Hesthaven1998}, acceptable values of $\tau$ span a broad range \cite{Hesthaven1998}.
   So it is possible to choose a $\tau$ that is suitable for all wavenumbers to make $B$ also independent of wavenumber.  
   Second, note that a division by $A - k_j^2I$ is required for each wavenumber $k_j$, and $A - k_j^2I$ certainly depends on wavenumber although only along its diagonal.
   To avoid operations of cubic complexity like Gaussian elimination to compute $(A - k_j^2)^{-1}$, or calculating and storing separate factorizations of $A(k_j)$ for all $j$, we appeal to the Schur factorization method (not to be confused with a Schur \emph{complement}).
   In particular, by computing and storing the Schur factorization of $A(0)$, divisions by $A(k_j)$ for all $k_j$ can be solved in quadratic time.
	Since the unshifted $A(0)$ is block diagonal, its Schur factorization $A(0) = UTU^*$ can be computed quickly in parallel.
   Recalling that in a Schur factorization $U$ is unitary and $T$ upper triangular, the Schur factorization of $A(k_j)$ for any $k_j$ is given by
	\begin{align}
		A(k_j) &= A -  k_j^2 I \\
		&= U T U^* - k_j^2 I \\
		&= U( T - k_j^2 I)U^*.
	\end{align}  
	$U$ is block unitary and $T - k_j^2I$ is block upper triangular, so any division $x = A(k_j)\backslash b$ is computable in parallel and in quadratic time by
	\begin{align}
		x = U( T - k_j^2 )^{-1} U^* b
	\end{align}
	where $(T - k_j^2 I)^{-1}$ is a block back-substitution and multiplications by $U$ and $U^*$ are block matrix-vector multiplications.  Both of these operations are parallel over the blocks of $A$ and are of quadratic complexity.  It is worth noting that although $A(k_j)\in\R^{r\times r}$, for unsymmetric matrices $U,T \in \C^{r \times r}$, which roughly doubles the storage requirement.

	From the perspective of domain decomposition, the Schur factorization is a powerful tool.
   It allows for solving linear systems in $A(k_j)$ in quadratic time while only requiring the Schur factorization and storage of $A(0)$.
   The Schur factorization of $A(0)$ is furthermore computable in parallel since $A(0)$, and consequently $U$ and $T$ are block-diagaonal.

   After computing $A(0) = UTU^*$, each transverse Schur problem is assembled as
   \begin{align}
      S(k_j) \hat x_j &= \hat b_j \\
      (I + BA(k_j)^{-1}E)\hat x_j \hat x_j &= \hat b_j \\
      (I + BU(T - k_j^2I)\backslash U^*E)\hat x_j &= \hat b_j,
   \end{align}
   where $\hat{x}_j = B\hat{u}_j$ and $\hat{b}_j = BA(k_j)^{-1}\hat{f}_j$. Then, each transverse wavenumber's solution $\hat u_j$ is obtained from the solution of the Schur problem $\hat{x}_j$ as
   \begin{align}
      \hat{u}_j &= A(k_j)^{-1}( \hat{f}_j - E\hat{x}_j) \\
      \hat{u}_j &= U(T - k_j^2I)^{-1}U^*( \hat{f}_j - E\hat{x}_j),
   \end{align}
   and the full three-dimensional solution is assembled via the inverse Fourier transform,
   \begin{align}
      u(x,y,z) = \sum_j \hat{u}_j(x,z) e^{ik_y y}.
   \end{align}
   Each Schur problem $S(k_j)\hat{x}_j = \hat{b}_j$ can be solved in parallel with the same deflation preconditioning technique demonstrated previously.
   The assembly of each Schur matrix $S(k_j)$ is non-trivial, but is amortized over the many Poisson solves required in a time-evolving Navier-Stokes solver and is in any case embarrassingly parallel over wavenumber.

	Thus, the deflated and preconditioned three-dimensional Schur problem is formulated as a sequence of independent two-dimensional Schur problems, each a   Helmholtz equation representing the $k_j$ wavenumber component of the solution $u(x,y,z)$.
   Identical to the formulation of the zero wavenumber two-dimensional problem in Eq.~(\ref{deflated_schur}), the deflated and preconditioned problems are for each $k_j$
	\begin{align}
		P(k_j)S(k_j)M(k_j) \hat x_j = P(k_j)BA(k_j)^{-1}\hat b_j
	\end{align}
	where $M(k_j)$ is the block-Jacobi preconditioner of $S(k_j)$, $P(k_j) = I + S(k_j)ZC(k_j)^{-1}Z^T$, and $C(k_j) = Z^T S(k_j)Z$.  The solution for wavenumber $k_j$ is obtained analogously to that in Eq.~(\ref{deflated_schur}) as
		\begin{align}
         \label{3d_deflated_schur}
			\hat{x_j} &= ZC(k_j) Z^T \hat{b}_j \nonumber \\ &+ Q(k_j)M(k_j)^{-1}\textrm{GMRES}(P(k_j)S(k_j)M(k_j)^{-1},P(k_j)\hat{b}_j).
		\end{align}
      where $Q(k_j) = I + ZC(k_j)^{-1}Z^T S(k_j)$.  The zero wavenumber problem $k_j = 0$ requires additional regularization as has been extensively documented in Section \ref{inconsistency}.

\subsection{Implementation}
   While each wavenumber Schur problem is independent of the others, in practice, all of the wave numbers are solved simultaneously in one GMRES calculation.
   This is done to minimize the communication overhead inherent in computation of dot products that plagues all Krylov methods \cite{Ghysels2013}; by solving all of the Schur problems together, the communication overhead penalty is paid only once per Krylov iteration instead of $m_y$-many times.

	The algorithmic summary of this method is given in Algorithm \ref{alg:3d}, which presumes that an initial setup phase has been conducted to compute the Schur factorization of $A(0)$, the assembly of $S(k_j)$ for $j = 1, \dots, m_y$, the block-Jacobi preconditioner $M(k_j)$, and the coarse grid matrix $C(k_j)$.  Since the two-dimensional problem corresponding to the $k_j = 0$ wavenumber is rank-deficient, a regularization identical to that in Algorithm \ref{alg:seq} is done for $k_y = 0$ (but is not shown in Algorithm \ref{alg:3d}).

	\renewcommand{\algorithmicrequire}{\textbf{Input:}}
	\renewcommand{\algorithmicensure}{\textbf{Output:}}
	\begin{algorithm}[h]
		\begin{algorithmic}[1]	
			\REQUIRE $b$
			\ENSURE $x$
			\FOR{$i = 1$ \TO $ n^2 m_xm_z $}
				\STATE $\hat b(i ,:) = $ FFT$[b(i,:)]$			
			\ENDFOR
			 \FOR{$j = 0$ \TO $ m_y/2 - 1$ } 			
			 	\STATE{ $\hat b(:,j) \leftarrow B U( T - k_j^2 I)^{-1}U^* \hat{b}(:, j )$ }
			\ENDFOR

				\STATE $\hat x = $ GMRES$( \sum_j P(k_j)S(k_j)M(k_j)^{-1},P \hat b)$
			\FOR{$j = 1$ \TO $ m_y$ } 			

				\STATE{$ \hat x(:,j) \leftarrow  U ( T - k_j^2 I)^{-1}U^*( b_S - E \hat{x}(:,j) )$}
			\ENDFOR			
			\FOR{$i = 1$ \TO $ n^2m_xm_z $}
				\STATE $x(i ,:) = $ IFFT$[\hat x(i,:)]$			
			\ENDFOR			
		\end{algorithmic}
		\caption{Deflated Schur complement method for the 3D Poisson equation.}
		\label{alg:3d}
	\end{algorithm}
	The two for-loops wrapping the fast Fourier transforms (FFTs) are embarrassingly parallel and in practice multi-threaded and computed with an external library\footnote{Furthermore, in the context of a Navier-Stokes simulation, the Poisson equation is solved entirely in Fourier space thus obviating the need for the Fourier transforms within the Poisson solver itself.}.  The divisions by $A(k_j)$ (lines 5 and 9) are block-solves each of quadratic complexity, and so are negligible relative to the expensive GMRES iteration for the solution of the three-dimensional Schur complement system (line 7).    When appropriately parallelized on $n_p$ processors, the number of floating point operations per rank for the above algorithm is
	\begin{align}
	\label{3dflops}
		F &= \overbrace{2n^2 m_z \frac{m_x} {n_p} \mathcal{O}(m_y \log m_y )}^{ \textrm{$n^2m_zm_x/n_p$ FFTs and IFFTs} } \nonumber \\ 
		&+ \underbrace{6 m_y\frac{m_x}{n_p}\mathcal{O}( (n^2 m_z) ^2 )}_{\textrm{ $6 m_y m_x/n_p$ Schur Solves}} + \underbrace{ \frac{m_x}{n_p} \mathcal{O}\left(  K^316 m_y n^2 m_z^2 \right)}_{ \textrm{1  GMRES solve of $K$ iterations}}
	\end{align}
	where $K$ is the number of GMRES iterations required for a solve.
   In practical applications $m_x /n_p$ is kept bounded, and so if $K$ depends on either $m_z$ or $n$, then the third term will quickly become dominant due to its cubic dependence on $K$.
   Thus, even in the three-dimensional case, minimizing GMRES iterations $K$ is still the most important aspect of achieving good performance.

\subsection{Performance}


	To study the performance of the three-dimensional solver, the two-dimensional simulations represented in Table \ref{table_time} were extended by extruding the two-dimensional domain into $\Omega \times [0,l_y]$, with $m_y$ grid points in the transverse direction.
   The right hand sides used are again randomly drawn from a uniform distribution to ensure that their Fourier transforms will have significant components in all wavenumbers $k_j$.
   As before, the number of elements in $x$ was increased from $m_x = 64$ to $m_x = 1024$, while the other parameters $(n,m_y,m_z) = (10, 32, 10)$ were held constant.
   These parameters imply grids ranging from 2 to 32 million grid points in size.
   From Table \ref{table_time} only the deflation and two-level additive Schwarz preconditioners were used, as the other two methods required too much computation time and Krylov iterations to be of practical use.
   All results were computed on $32$ processors, and all results were averaged over $10$ trials.


	The results comparing the deflation preconditioning with two-level additive Schwarz preconditioning are displayed in Table \ref{table_3dtime}.   First notice that the number of GMRES iterations does not increase substantially from the two-dimensional problems (Table \ref{table_time}).  This is primarily because the worst-conditioned of all the wavenumbers is the $k_0 = 0$ wavenumber which is what is solved in the two-dimensional problems; adding more wavenumbers only adds better conditioned problems that converge earlier than the zero wavenumber.  The second thing to note is that the deflation method again outperforms two-level additive Schwarz preconditioning, requiring about $40\%$ fewer GMRES iterations and $15\%$ less time.  The difference in iterations does not translate directly into time because again since the application of $P$ in the deflation method adds an extra multiplication of $S$ to each GMRES iteration.  While, from one perspective, deflation is only modestly faster than two-level preconditioning, when solving the Poisson equation hundreds if not thousands of times within a time-evolving fluid simulation, even a $15\%$ improvement in time translates to significant savings in total time.  Besides which, for extremely large problems for which storage of Krlyov basis vectors is a significant cost the reduction in GMRES iterations in deflation relative to two-level additive Schwarz preconditioning can be important.
	
\begin{table}[h]
   	\begin{center}
	\caption{Comparison between the 3D deflation method and the 3D two-level additive Schwarz method of iterations and computation time to solution in seconds.  The number of GLL points per direction is $n = 10$, the number of vertical elements $m_z = 10$, and number of transverse wave numbers $m_y =32$.  These simulations were all benchmarked in an MPI-parallel Fortran code and executed on 32 processors.}
	\label{table_3dtime}
	\begin{small}
	\scalebox{0.8}{
        \begin{tabular}{ r  r  c | c  c | c  c |  }        
            \cline{4-7}
             & & & \multicolumn{2}{|c|}{Deflation}& \multicolumn{2}{|c|}{TL Schwarz}\\
            $m_x$  & Grid Points & Setup Time & Iter. & Time & Iter. & Time \\
            \hline
            \multicolumn{1}{|r|}{64}     &  $2.05 \times 10^7$    & 7.50e1 &  30.5    & 2.30e0    & 48.1    & 2.71e0   \\
            \multicolumn{1}{|r|}{128}   &  $4.10 \times 10^7$    & 1.52e2 &  57.0     &  7.21e0  & 84.1   & 8.29e0 \\
            \multicolumn{1}{|r|}{256}   &  $8.19 \times 10^7$    & 3.13e2 &  39.2     & 1.08e1    & 60.0  & 1.23e1 \\
            \multicolumn{1}{|r|}{512}   &  $1.64 \times 10^8$    & 6.58e2 &  33.7     & 1.98e1   & 52.7   & 2.27e1 \\
            \multicolumn{1}{|r|}{1024} &  $3.28 \times 10^8$    & 1.50e3 &  31.8      & 5.21e1   & 50.4   & 6.50e1 \\
            \hline
        \end{tabular}
        }
        \end{small}
        \end{center}
  \end{table}		
  
  Lastly note that the third column in Table \ref{table_3dtime} shows the time taken to perform the assembly, factorization, and inverse iteration done in the setup phase initially.  The time for setup is usually about 20 times the time taken for a single solve.  This means that so long as the Poisson problem is to be solved many more than 20 times, the cost of setup is negligible.  In applications in environmental flows, there are usually $\mathcal{O}(10^4-10^5)$ many time-steps each requiring the solution of a Poisson-Neumann problem which amortizes the cost of setup over enough solves to make the setup cost acceptable \cite{Diamessis2005, Escobar-Vargas2014}.
  
\section{Discussion}

The SMPM was used here as the platform to demonstrate the effectiveness of the deflated/preconditioned Schur complement approach, but this approach can be extended to other discretizations as well.
Domain decomposition has of course been applied to many different kinds of discretizations\cite{Toselli2004} including continuous \cite{Pavarino2000} and discontinuous Galerkin, collocation \cite{Bialecki2007}, and spectral element methods \cite{Manna2004, Couzy1995}.
While most domain decomposition methods separate the grid into internal and interface unknowns, the particular decomposition shown in Section \ref{schur_assembly} decomposes the discrete operator into intra- and inter-subdomain components and solves for the inter-subdomain fluxes $Bu$ first before solving the local block problems.
This approach is similar to one taken in a class of methods known as Finite Element Tearing and Interconnect (FETI) \cite{farhat1991}.
In FETI, a set of Lagrange multipliers that represent inter-subdomain fluxes are solved first, and these Lagrange multipliers are the direct analogues of the inter-subdomain fluxes $Bu$ in SMPM.
Unlike in FETI, the element matrices in the SMPM are all invertible (in FETI elements without no intersection with the outer boundaries have non-invertible element matrices) and the decomposition here does not require that the operator be symmetric positive definite as FETI requires for the construction of its weak formulation.
Indeed the SMPM operator matrices are neither positive definite nor symmetric, and there is no weak form owing to the collocation-based nature of the SMPM.

While the domain decomposition, construction of the Schur problem deflation vectors, and block-Jacobi preconditioner used here are generally applicable to any element-based discretization of an elliptic problem, some properties inherent to the SMPM are required.
First, it is critical that the discretization be discontinuous, as the Schur problem constructed in Section \ref{schur_assembly} assumes that the solution can be discontinuous across the decomposed interfaces.
Because the SMPM does not invoke strong continuity across elements the operator $A$ is block-diagonal and hence the interior problems (Step 2 and Step 9 in Algorithm \ref{alg:seq}) can be solved in parallel.
Lastly, and perhaps most importantly, because the SMPM is discontinuous the local block problems in $A$ are each well-defined, representing totally decoupled homogenous Robin boundary value problems; if the SMPM were continuous the decomposition would have been less clearly defined.

Secondly, it is helpful that the SMPM when used in practice is of high-order ($n\gg 1$), which makes the Schur problem significantly smaller than the full problem.
Since the dimension of $L$ grows $O(n^2)$ with $n$ and the dimension of $S$ grows as $O(n)$, the savings in storage of the Krylov basis as well as in the floating point operations in orthogonalization can be of importance.
In practice, values as high as $n=32$ have been used\cite{Diamessis2006}, but more typically $n \leq 20$.

These two tenants (discontinuous and high-order) are satisfied by the high-order discontinuous Galerkin (DG) class of methods\cite{Cockburn2000}, and using the deflated/preconditioned Schur complement approach described in this work seems practicable.  The boundary fluxes $Bu$ in the SMPM that are the unknowns of the Schur problem correspond to numerical fluxes in DG, but otherwise the decomposition and solution would proceed identically.
Since DG methods are widely used, this seems like a natural application for deflation-augmenting preconditioning of the Schur complement problem.

Finally, while this paper focused solely on the Poisson-Neumann problem due to its relevance in the problem of interest, the approach taken extends to any elliptic problem.
Thus the only stringent requirements are that the partial differential equation be elliptic, and the method be element-based, discontinuous, and preferably of high-order.

\section{Summary and Concluding Remarks}

In this paper, a method of solving the Poisson-Neumann problem on large, highly-elongated domains by way of a domain decomposition is presented.
In particular, a method of preconditioning the Schur complement problem of the two-dimensional Poisson-Neumann system discretized by a high-order discontinuous spectral element method has been developed.
The large, highly-elongated grids introduce two difficulties to solving the Poisson-Neumann system.
First, the highly-stretched grids lead to ill-conditioned spectral differentiation matrices; this was addressed by way of a block-Jacobi preconditioner on the Schur problem.
Second, for grids with many elements, Krylov subspace methods require many iterations to communicate information across the grid; this difficulty was addressed by augmenting the block-Jacobi preconditioner with a set of deflation vectors that eliminated the residual in a coarse subspace of the Schur matrix.
Since the Schur complement problem is better conditioned and a factor of $n$ smaller than the full Poisson problem (where $n+1$ is the order of the polynomial basis in each element), it requires fewer Krylov iterations and less storage for the Krylov basis vectors.
It was shown that by using deflation-augmented preconditioning on the Schur complement, the number of Krylov iterations can be kept bounded as the aspect ratio was increased and as the number of $x$-elements is increased.
Comparisons with two-level additive Schwarz preconditioning for the same problem showed that deflation required approximately half as many Krylov iterations and required about $25\%$ less time.
Lastly it was shown that when a third dimension was added by way of a Fourier discretization, the superiority of deflation in two-dimensions translated to the three-dimensional problem as well; in three dimensions, the deflation preconditioned Schur problem required about $40\%$ fewer Krylov iterations and about $15\%$ less time than two-level additive Schwarz preconditioning.



The three main contributions of this work may be summarized as follows.
First, the use of deflation-augmented block-Jacobi preconditioning on the Poisson-Neumann Schur system is shown to be a viable candidate technique for obtaining the solution to the Schur complement problem.
This method is viable in the sense that convergence of the GMRES calculation on the deflated and preconditioned Schur matrix is independent of the domain aspect ratio and the number of $x$ elements, both important criteria for a solution algorithm and the subject of much study \cite{Fischer1997,Fischer2005,Feng2001}.
Development of algorithms that satisfy these criteria has largely been focused on two-level additive Schwarz methods applied to the Poisson problem itself.
While there admittedly exists a significant body of literature on using two-level additive Schwarz preconditioning for the Schur problem of elliptic equations \cite{Pavarino2000,Bialecki2007,Manna2004}, and also a significant body of literature on using deflation methods to accelerate convergence on elliptic problems themselves \cite{Vuik2006,Nabben2004,Tang2012}, applying deflation acceleration to the Schur complement problem seems largely unexplored.
It is also worth noting that the SMPM discretized Poisson-Neumann operator is unsymmetric, and some care was taken to show that the projections out of the null space of this rank-deficient operator did not affect the quality of the solution.
While certainly most discretizations of the self-adjoint Poisson equation are symmetric positive definite themselves (as in the spectral element or discontinuous Galerkin methods), this also shows that, without much extra effort, unsymmetric discretizations can also be accounted for.

Secondly, it is demonstrated in Sections \ref{sec:2d} and \ref{sec:3d} that deflation acceleration alongside block-Jacobi preconditioning outperforms two-level additive Schwarz at least for the SMPM discretization of the Poisson-Neumann problem.
For both the two-dimensional (Section \ref{sec:2d}) and the three-dimensional (Section \ref{sec:3d}) case, deflation/block-Jacobi preconditioning required about half as many GMRES iterations and between 15\% and 25\% less computation time than two-level additive Schwarz.
While this is admittedly a limited comparison and a modest performance gain, it is in agreement with previous work comparing deflation with two-level additive Schwarz methods for preconditioning the full Poisson problem (i.e. not the Schur complement) \cite{Nabben2004, Tang2010, Tang2009}.
And in light of the fact that for many of these environmental-scale problems the duration of the computation is on the order of several months of wall-clock time \cite{Diamessis2011}, even a $15\%$ reduction in computation time as shown in Section \ref{sec:3d} can be significant.
Considering also the $50\%$ reduction in the memory footprint due to the $50\%$ reduction in Krylov iterations (and thus basis vectors), it is possible that the performance difference between the deflation and the two-level Schwarz methods may be significant.

Finally, use of the Schur factorization of the local block problems $A(k_j)$ is used to quickly factor and divide all the local element matrices in transverse direction when extending to three-dimensions (Section \ref{sec:3d}).
Coupled with the domain decomposition outlined in Section \ref{sec:2d}, the use of the Schur factorizations allows for parallel division of the block-diagonal matrix $A(k_j)$ and the re-use of the Schur factorization of $A(0) = UTU^*$. 
Without the Schur factorization, either the factorization of each $A(k_j)$ would need to be stored or a different approach such as Gaussian elimination (which is of cubic complexity) would be required to perform the divisions of $\{A(k_j)\}_{j=1}^{m_y}$.
For grids with many transverse wavenumbers ($m_y \gg 1$) the savings in storage and factorization time can be immense, and in practice $m_y$ is often much larger than the $m_y = 32$ value used in Section \ref{sec:3d}, reaching from $m_y = 128$ \cite{Diamessis2005} to as high as $m_y = 512$ \cite{Zhou2013} in supercomputing-scale applications.
When $m_y$ is this large, storing a single $A(0)$ and computing a single Schur factorization instead of $m_y$-many is an important advantage of a Fourier-discretized transverse direction leveraging the Schur factorization.

\subsection{Future work}

Since the primary motivation for solving the Poisson equation in this work is to obtain a numerical solution to the incompressible Navier-Stokes equations,  immediate interest is in implementing the deflated block-Jacobi Schur method in an incompressible Navier-Stokes solver.
By way of the operator splitting as in Ref.~\cite{Karniadakis1991}, the Poisson problem is separated from the viscous and nonlinear advective terms which means the methodology developed in this work can be applied in a relatively straightforward manner.

Some extensions to the deflation block-Jacobi method described here might be to a fully three-dimensional (i.e. in which elements are mapped from $[-1,1]^3$) discretization or to an unstructured grid.
First, in a fully three-dimensional grid, the Schur complement would likely be too large too store, but otherwise the construction of its deflation vectors and its preconditioner would be an extension of the current aproach.
For example if the Schur matrix was constructed by decomposing along all element interfaces, each block of the local matrix $A$ would be of dimension $n^3$, and the Schur complement matrix $S$ would be of dimension $6n^2(m_x - 1)(m_y-1)(m_z-1)$.  
Thus the Schur matrix would still be a factor of $n$ smaller than the full Poisson matrix, and would be better conditioned by virtue of being a Schur complement.
If in addition to being fully three-dimensional the grid were also unstructured the choice of subdomains becomes more complicated, but can be reduced to a graph partitioning problem that has received much attention as a computational problem \cite{Karypis2013}.
Thus, while requiring a significantly greater programming investment, the extension of the Schur complement and deflation methods outlined in this paper to an unstructured three-dimensional grid is theoretically tractable.

One open question worth investigating is related to the choice of deflation vectors $Z$ in Section \ref{sec:def}.  
The deflation vectors were chosen to be discrete indicator vectors along the domain decomposition boundaries $\Gamma_i$ after a choice made in Section 4.1.1 of Ref.~\cite{Tang2009}, a simulation of porous media flow that required the solution of Poisson-Dirichlet system of equations.
This choice is sensible, as it approximates eigenvectors that represent modes of the Schur matrix in which all nodes along an interface $\Gamma_i$ have the same value; this represents low-order modes of oscillation in a sense, and is in the spirit of what deflation is attempting to do.
However, we do not claim (or even suggest) that this choice of $Z$ is optimal, and it is likely that other deflation vectors and thus coarsening operators will provide better performance.

Lastly, although it is true that by replacing the boundary fluxes $Bu$ with the inter-element fluxes in the discontinuous Galerkin method (DGM) the Schur complement assembly and deflation above can extend to high-order DGM discretizations, a demonstration of that claim with the aim of studying the relative performance of the SMPM against the DGM would be highly insightful.

\section*{Acknowledgements}

The authors would like acknowledge Professors Charles Van Loan, Jorge Escobar-Vargas, and David Bindel for helpful comments on the work in this manuscript.  The authors would also like to thank the United States Department of Defense High-Performance Computing Modernization Office for the NDSEG fellowship and the National Science Foundation for CAREER award \#0845558 for support.

\section*{Appendix: proof of claims}

   \subsection*{Proof of Claim 1}
	\begin{proof}
		First, note that $u_L^T(A + EB) = 0_r \in \R^{r}$ which implies that $u_L^T + u_L^T EBA^{-1} = 0_r$.  Right multiply by $E$ to obtain $u_L^TE + u_L^TEBA^{-1}E = 0_r$ which yields the first relation that
		$u_L^TE( I + BA^{-1}E) = 0_r \Rightarrow u_L^TES = 0_k \Rightarrow u_S = E^Tu_L \in \R^k$.

		For the second part, note that $BA^{-1}L = BA^{-1}(A + EB) = B + BA^{-1}EB = SB$.  Thus if $u_S^TS = 0_k$, $u_S^TBA^{-1}L = 0_k$, and so $A^{-T}B^Tu_S \in \R^r$ is the left null vector of $L$.
	\end{proof}

   \subsection*{Proof of Claim 2}
	\begin{proof}
		Start with the Poisson residual $Lu - \tilde{f}$, substitute the solution $u = A^{-1}(\tilde{f} - Ex_s)$, and recall that $L = A + EB$ to obtain
		\begin{align}
			Lu - \tilde{f} &= LA^{-1}(b - Ex_s) - \tilde{f} \\
					  &= (A + EB)A^{-1}(\tilde{f} - Ex_S) -\tilde{f} \\
					  &= (I + EBA^{-1})(\tilde{f} - Ex_S) -\tilde{f} \\
					  &= EBA^{-1}\tilde{f} - ESx_S \\
					  &= E( BA^{-1}\tilde{f} - Sx_S).
		\end{align}
		If $S$ were full-rank, then $x_s = S^{-1}BA^{-1}\tilde{f}$ as calculated in domain decomposition and the right hand side would be zero.  However due to the rank-deficiency of $S$ the residual that is made small in the Schur complement solve is $\norm{ Sx_s - (I - u_Su_S^T) BA^{-1}\tilde{f} }_2$.  So, adding and subtracting the projection to the above, we obtain a bound on the residual:
		\begin{align}
		\norm{Lu - \tilde{f}}_2&= \norm{E(BA^{-1}\tilde{f} - u_Su_S^T BA^{-1}\tilde{f} + u_Su_S^TBA^{-1}\tilde{f} - Sx_s)}_2 \\
		&\le \norm{E(BA^{-1}\tilde{f} - u_Su_S^T BA^{-1}\tilde{f} - Sx_s)}_2 + \norm{Eu_Su_S^TBA^{-1}\tilde{f}}_2 \\
		&\le \norm{ (I - u_Su_S^T)BA^{-1}\tilde{f} - Sx_s)}_2 + \norm{Eu_Su_S^TBA^{-1}\tilde{f}}_2.
		\end{align}
		The first term in the above is exactly the residual we minimize in the Schur complement solve; the second term is an error term related to projecting out of the column space of $S$.  Examining the second term, recall from Claim 1 that $u_S = E^T u_L$, and substitute this into the expression for the second term:
		\begin{align}
			\norm{Eu_Su_S^TBA^{-1}\tilde{f}}_2 &= \norm{ E (E^Tu_L)(E^T(u_L))^T BA^{-1}\tilde{f} }_2 \\
			&= \norm{ E (E^Tu_L)u_L^T E BA^{-1}\tilde{f} }_2.
		\end{align}
		Recalling that $EB = L - A$ we get
		\begin{align}
			&= \norm{ E E^Tu_L u_L^T( L - A )A^{-1}\tilde{f} }_2 \\
			&=  \norm{ E E^Tu_L ( u_L^T LA^{-1}\tilde{f} - u_L^T \tilde{f} ) }_2.
		\end{align}
		Now by definition $u_L^TL = 0_r$ so the first term in the above is zero, and by assumption $u_L^T\tilde{f} = 0_r$ so the second term is zero.  Thus $\norm{ Eu_Su_S^TBA^{-1}\tilde{f} }_2 = 0$, and we have our bound,
		\begin{align}
			\norm{ Lu - \tilde{f}}_2 \leq \norm{Sx_S - (I - u_S u_S^T)b_S }_2,
		\end{align}
		where $b_S = BA^{-1}\tilde{f}$.
	\end{proof}

\bibliographystyle{elsarticle-harv}
\bibliography{capacitance_bibliography}

\end{document}